\documentclass{gtart}
\gtart

\usepackage{amsthm}
\usepackage{amsgen}
\usepackage{amsmath,amssymb}
\usepackage{curvesls}
\usepackage{epic}
\usepackage[all]{xy}
\usepackage[dvips]{graphicx}
\usepackage{psfrag}
\usepackage[T1]{fontenc}
\usepackage[sc]{mathpazo}
\usepackage{color, xcolor, colortbl}
\usepackage{float}
\usepackage{arydshln, multirow} 
\usepackage{multicol}
\usepackage{fancybox}
\newtheorem{thm}{Theorem}[section]
\newtheorem{lem}[thm]{Lemma}

\newtheorem{prop}[thm]{Proposition}

\theoremstyle{definition}
\newtheorem{defn}[thm]{Definition}

\newcommand{\Real}{{\mathbb R }}
\newcommand{\Rat}{{\mathbb Q}}
\newcommand{\Zed}{{\mathbb Z }}

\newcommand{\Emb}{{\mathrm{Emb}}}

\newcommand{\Diff}{{\mathrm{Diff}}}
\newcommand{\Homeo}{{\mathrm{Homeo}}}
\newcommand{\PDiff}{{\mathrm{PDiff}}}

\newcommand{\Map}{{\mathrm{Map}}}

\newcommand{\fp}{{\mathcal{F}}}

\DeclareMathOperator{\id}{id}

\newtheorem{proposition}[thm]{Proposition}

\newcommand{\BZ}{\mathbb Z}
\newcommand{\BH}{\mathbb H}

\definecolor{RedCol}{rgb}{1.0,0.5,0.5}
\definecolor{BlueRow}{rgb}{0.7,0.7,1.0}
\newcolumntype{R}{>{\columncolor{RedCol}}c}

\definecolor{orange}{RGB}{255, 123, 21}

\begin{document}

\title{On the automorphism groups of hyperbolic manifolds}
\authors{Ryan Budney \\ David Gabai}

\addresses{
Mathematics and Statistics, University of Victoria PO BOX 3060 STN CSC, Victoria BC Canada V8W 3R4\\
Fine Hall, Washington Road Princeton NJ 08544-1000 USA}
\emails{rybu@uvic.ca \\ gabai@math.princeton.edu}

\begin{abstract} 
Let $\Diff(N)$ and $\Homeo(N)$ denote the smooth and topological group of automorphisms respectively that
fix the boundary of the $n$-manifold $N$, pointwise.  We show that 
$\pi_{n-4} \Homeo(S^1 \times D^{n-1})$ is not finitely-generated
for $n \geq 4$ and in particular $\pi_0\Homeo(S^1\times D^3)$ is infinitely generated.  
We apply this to show that the smooth and topological automorphism groups of finite-volume
 hyperbolic $n$-manifolds (when $n \geq 4$) do not have the homotopy-type of finite CW-complexes, 
 results previously known for $n\ge 11$ by Farrell and Jones.  In particular, we show that if $N$ 
 is a closed hyperbolic $n$-manifold, and $\Diff_0(N)$ represents the subgroup of diffeomorphisms
 that are homotopic to the identity, then $\pi_{n-4}\Diff_0(N)$ is infinitely generated and 
 hence if $n=4$, then $\pi_0\Diff_0(N)$ is infinitely generated with similar results holding 
 topologically.  \end{abstract}

\primaryclass{57M99}
\secondaryclass{57R52, 57R50, 57N50}
\keywords{4-manifolds, 2-knots, isotopy}
\maketitle


\section{Introduction}\label{intro}

The main result of this paper is the following.

\begin{thm} \label{main} $\pi_{n-4} \Homeo(S^1 \times D^{n-1})$ is infinitely generated and in particular $\pi_0\Homeo(S^1\times D^3)$ is infinitely generated. \end{thm}

  In the smooth category, this was the topic of \cite{BG}, where it was shown that
$\pi_{n-4} \Diff(S^1 \times D^{n-1})$ is not finitely generated.  Here all automorphism groups  act via the identity on the boundary and hence a given automorphism is homotopic to $\id$.  To prove this theorem we elaborate a method briefly introduced in \cite{BG} using linking numbers coming from collinear and cohorizontal spaces and use it to give in Section \ref{deltk} a new  proof that the $\delta_k$ families of \cite{BG} are linearly independent in the smooth category, provided $k\ge 4$.    It's new in the sense that it is a direct argument using that theory starting with the $\delta_k$ families while \cite{BG} showed how to express $\delta_k$'s in terms of our $G(p,q)$ families.  Remarkably, being based on elementary intersection theory, this method also works in the topological category  as detailed in Section \ref{homeosec}.  Our main result has the following applications.

\begin{thm}  The automorphism groups of $S^1 \times D^{n-1}$ do not have the homotopy-type of finite-dimensional 
CW-complexes, provided $n \geq 4$.\end{thm}

For dimensions $n\ge 6 $ this result was  proven by Hatcher and Wagoner \cite{HW} more than 50 years ago, where they showed that the topological and smooth mapping class groups of $S^1 \times D^{n-1}$
are not finitely generated.  In contrast, the smooth and topological automorphism groups of $S^1 \times D^{n-1}$ have the
homotopy-type of $\Omega S^1 \simeq \Zed$ when $n=2$, and when $n=3$ these groups are contractible by the work of Hatcher \cite{H2}.

\begin{thm}\label{hyperbolic}  If $N$ is a complete hyperbolic $n$-manifold, then $\pi_{n-4} \Diff_0(N)$ and $\pi_{n-4}\Homeo_0(N)$ are infinitely generated.  In particular if $n=4$, both $\pi_0\Diff_0(N)$ and $\pi_0\Homeo_0(N)$ are infinitely generated.\end{thm}

Here $\Diff_0$ and $\Homeo_0$ denote automorphisms homotopic to $\id$. For $n\ge 11$, this result was  proven by Farrell and Jones \cite{FaJo} over 30 years ago.   Our result is sharp since $\Diff_0(N)$ is contractible when $n\le 3$ by \cite{Ga2} and \cite{Gr}.  Details are given in Section \ref{hypsec}

Section \ref{bbdif} introduces the barbell manifolds $\mathcal B_{i,j}^n$ and defines corresponding barbell diffeomorphisms generalizing the notion of barbell manifolds and  diffeomorphisms given in \cite{BG}.  They are of independent interest and will be explored in a future paper.  The barbell manifolds of \cite{BG} are the ones denoted here by $\mathcal B_{n-2, n-2}^n$.  This section includes two definitions of barbell diffeomorphisms, one from a perspective analogous to the definition of a Dehn twist, in terms of resolutions of double-points. The other perspective uses a product decomposition of barbell manifolds, and constructs barbell diffeomorphisms as fibre-preserving
maps, classically known as `horizontal' diffeomorphisms.  We also include some constructions of null isotopies, and null pseudoisotopies, 
for certain implanted barbell diffeomorphisms.

\noindent\emph{Acknowledgements}: Part of this work was developed during two visits to BIRS and one to Oberwolfach.  We thank these institutions for their hospitality.  The authors would like to thank Allen Hatcher, Sander Kupers and Manuel Krannich for helpful comments on an initial
draft of this paper. 


\section{Barbell diffeomorphisms}\label{bbdif}

For the purpose of this paper, an $n$-dimensional barbell manifold will be the boundary connect-sum of two trivial disc-bundles over spheres. 
We index the barbell manifolds by the dimensions of the spheres, thus we define
$$\mathcal B_{i,j}^n = S^i \times D^{n-i} \natural S^j \times D^{n-j}$$
as the standard $(i,j)$-barbell in dimension $n$.  We will always assume $i,j \geq 1$, as none of our constructions below
will be of interest when $i=0$ or $j=0$.  The spheres $S^i \times \{0\}$ in the first summand and $S^j \times \{0\}$ in the
second summand we call {\it core spheres}.  The discs $\{*\} \times D^{n-i}$ in the first summand and $\{*\} \times D^{n-j}$
in the second we call the {\it cocores}, where $\{*\}$ is a choice of basepoint in the respective spheres.  The {\it mid-ball}
we denote $B^{n-1}$, this is the embedded co-dimension one disc that separates the boundary connect sum into a copy of
$S^i \times D^{n-i}$ and $S^j \times D^{n-j}$ respectively. 

We will use the terminology $\Diff(M)$ to denote the group of diffeomorphisms of a manifold.  If $M$ has
boundary, we demand the diffeomorphisms restrict to the identity on the boundary, i.e. the restriction map
$\Diff(M) \to \Diff(\partial M)$ is a constant function. 

For the sake of argument, assume $i \leq j$. Consider the barbell manifold as fibering over $D^{n-j-1}$ with 
fiber $\mathcal B_{i,j}^{j+1}$.  As a manifold, $\mathcal B_{i,j}^{j+1}$ is the once-punctured $S^i \times D^{j-i+1}$.  
If we let $\Diff^\fp(\mathcal B_{i,j}^n)$ denote the fiber-preserving diffeomorphism group of $\mathcal B_{i,j}^n$, i.e,
diffeomorphisms $f : \mathcal B_{i,j}^n \to \mathcal B_{i,j}^n$ giving commutative diagrams

$$\xymatrix{\mathcal B_{i,j}^n \ar[rr]^f \ar[dr] && \mathcal B_{i,j}^n \ar[dl] \\
 & D^{n-j-1} & }.$$
 
These are sometimes also known as `horizontal diffeomorphisms.'
 Thus we have a homotopy-equivalence 
 $\Diff^\fp(\mathcal B_{i,j}^n) \simeq \Omega^{n-j-1} \Diff(\mathcal B_{i,j}^{j+1})$. 

Given that $\mathcal B_{i,j}^{j+1}$ is a once-punctured $S^i \times D^{j-i+1}$, there is the restriction fibre-bundle 
$$\Diff(\mathcal B_{i,j}^{j+1}) \to \Diff(S^i \times D^{j-i+1}) \to \Emb(D^{j+1}, S^i \times D^{j-i+1})$$
where the map to the base space is null-homotopic.  The map is obtained by fixing a compact $(j+1)$-ball in the
interior of $S^i \times D^{j-i+1}$ and taking the restriction map from $\Diff(S^i \times D^{j-i+1})$. 
Thus we have a fibre sequence 
$$\Omega \Emb(D^{j+1}, S^i \times D^{j-i+1}) \to \Diff(\mathcal B_{i,j}^{j+1}) \to \Diff(S^i \times D^{j-i+1})$$
such that the induced maps on homotopy groups give short exact sequences
$$0 \to \pi_k \Omega \Emb(D^{j+1}, S^i \times D^{j-i+1}) \to \pi_k \Diff(\mathcal B_{i,j}^{j+1}) \to \pi_k \Diff(S^i \times D^{j-i+1}) \to 0.$$
The map $\Omega \Emb(D^{j+1}, S^i \times D^{j-i+1}) \to \Diff(\mathcal B_{i,j}^{j+1})$ is obtained by applying isotopy extension
to the loop in $\Emb(D^{j+1}, S^i \times D^{j-i+1})$, and restricting to $\mathcal B_{i,j}^{j+1}$, i.e. the punctured copy of
$S^i \times D^{j-i+1}$. 

By the (homotopy) classification of spaces of tubular neighbourhoods, we have that $\Emb(D^{j+1}, S^i \times D^{j-i+1})$
has the homotopy-type of $S^i \times O_{j+1}$.  Observe that the generator of $\pi_i S^i \simeq \Zed$ gives a non-torsion
element of $\pi_{i-1} \Diff(\mathcal B_{i,j}^{j+1})$ via the above short exact sequence with $k=i-1$, and
via the equivalence $\Diff^\fp(\mathcal B_{i,j}^n) \simeq \Omega^{n-j-1} \Diff(\mathcal B_{i,j}^{j+1})$ it gives us 
a non-torsion element in $\pi_{i+j-n} \Diff^\fp(\mathcal B_{i,j}^n)$ provided $i+j \geq n$. 

We now analyze three special classes which are not completely disjoint:

\begin{enumerate}
\item In the case of the twice punctured $2$-disc $\mathcal B_{1,1}^2$, the barbell diffeomorphism is the composite
of the Dehn twists \cite{Dehn} about the boundary curves such that the signs form a homology, i.e. signs chosen consistent with
the boundary orientation.
\item The barbell diffeomorphism of $\mathcal B_{n-2,n-2}^n$ is the family studied in \cite{BG}. These barbells have the feature
that one can knot them in the `handcuff' fashion, provided $n \geq 3$. The diffeomorphisms
themselves are defined only when $i+j \geq n$, thus requires $n \geq 4$. 
\item When $i+j=n$ these barbells can be `Hopf-linked' in $S^1 \times D^{n-1}$ 
provided $i,j \geq 3$, i.e. $n \geq 6$, allowing us to relate to the work of Hatcher and Wagoner \cite{HW}. 
\end{enumerate}


We offer an alternative, more symmetric definition of the induced map 
$\pi_{i+j-n} \Omega^{n-j} S^i \equiv \Zed \to \pi_{i+j-n} \Diff(\mathcal B_{i,j}^n)$ when $i+j \geq n$. 
Consider two vector subspaces of $\Real^n$ isomorphic to $\Real^i$ and $\Real^j$.  We assume the two vector subspaces meet
in a single point, $\{0\}$, thus $n \geq i+j$.  If $n > i+j$ we can use a small perturbation near the origin (say, using 
bump function) to deform the vector subspaces to disjoint submanifolds.  Provided $n > i+j+1$, all such small deformations
are isotopic, as the normal sphere to the subspace spanned by $\Real^i$ and $\Real^j$ is $S^{n-i-j-1}$, which is 
connected. To ensure we are dealing with compact manifolds, consider $D^i \subset \Real^i$ and $D^j \subset \Real^j$. Using
bump functions supported in the interiors of these discs, gives us the following proposition.

\begin{figure}[H]
{
\psfrag{ri}[tl][tl][0.7][0]{$\Real^i$}
\psfrag{rj}[tl][tl][0.7][0]{$\Real^j$}
\psfrag{c}[tl][tl][0.7][0]{$\Real^{n-i-j}$}
$$\includegraphics[width=10cm]{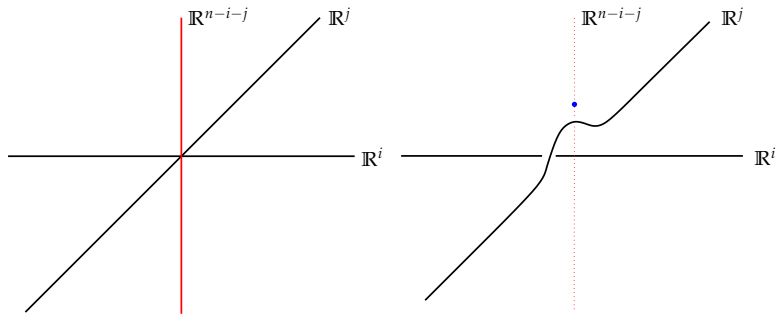}$$
\caption{Barbell diffeomorphism via resolution of double point.}\label{bbdoublept}
}
\end{figure}

\begin{prop}\label{twobb} Consider the spherical family $S^{n-i-j-1} \to \Emb(D^i \sqcup D^j, D^n)$ defined above. Then the 
connecting map for the 
homotopy long exact sequence for the fibration $\Diff(D^n) \to \Emb(D^i \sqcup D^j, D^n)$ gives us an element of
$\pi_{n-i-j-2} \Diff(D^n, D^i \sqcup D^j)$.  By thickening the embedded copies of $D^i$ and $D^j$ slightly, we can assume
these diffeomorphisms are the identity in a neighbourhood of the embedded copies of $D^i$ and $D^j$, thus this is an element
of the homotopy group
$$\pi_{n-i-j-2} \Diff(\mathcal{B}_{n-i-1, n-j-1}^n).$$
Moreover, if we let $i'=n-i-1$ and $j'=n-j-1$, this can be rewritten as an element of 
$\pi_{i'+j'-n} \Diff(\mathcal{B}_{i',j'}^n)$, and it is the barbell diffeomorphism, i.e. the induced map
on $\pi_{i'+j'-n}$ for the map $\Omega^{n-j'} S^{i'} \to \Diff(\mathcal B_{i',j'}^n)$. 
\end{prop}

String link families in Proposition \ref{twobb} are studied systematically in Koytcheff \cite{K}.


\begin{prop}\label{bbtriv}
The barbell diffeomorphism $\Zed \equiv \pi_{i+j-n} \Omega^{n-j} S^i \to \pi_{i+j-n} \Diff(\mathcal{B}_{i,j}^n)$
is essential, i.e. there is a homomorphism
$$\pi_{i+j-n} \Diff(\mathcal{B}_{i,j}^n) \to \Zed$$ 
that detects the barbell diffeomorphism.  This homomorphism is a version of the {\bf scanning map}.  Specifically, 
let $B$ be a mid-ball for $\mathcal{B}_{i,j}$, i.e. a smoothly-embedded copy of $D^{n-1}$ that
splits $\mathcal{B}_{i,j}^n$ into a boundary connect-sum.  Fiber $B$ by parallel intervals.  Scanning
using $B$ gives a map
$$\pi_{i+j-n} \Diff(\mathcal{B}_{i,j}^n) \to \pi_{i+j-2} \Emb(I, \mathcal{B}_{i,j}^n).$$
This map detects the barbell diffeomorphism.  Furthermore, the homomorphism
$$\pi_{i+j-2} \Emb(I, \mathcal{B}_{i,j}^n) \to \Zed$$
is computed by counting signed pairs of points $t_1 < t_2 \in I$ such that $f(t_1)$ is on the first cocore, and 
$f(t_2)$ is on the second cocore.  
\end{prop}

\begin{figure}[H]
{
\psfrag{SUM}[tl][tl][1][0]{$\mathcal{B}_{n-2,n-2}^n$}
\psfrag{E1}[tl][tl][0.7][0]{\textcolor{red}{$E_1$}}
\psfrag{E2}[tl][tl][0.7][0]{\textcolor{blue}{$E_2$}}
\psfrag{DN}[tl][tl][0.8][0]{$D^{2n-6} \equiv D^{n-3} \times D^{n-3}$}
\psfrag{Dn}[tl][tl][0.8][0]{$D^{n-3}$}

$$\includegraphics[width=10cm]{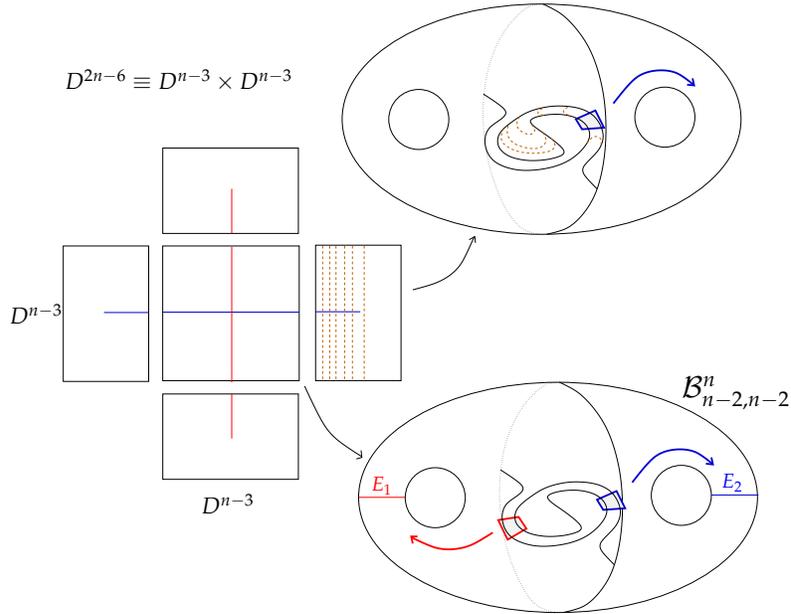}$$
\caption{Barbell diffeomorphism family restricted to mid-ball as
map $D^{2n-6} \to \Emb(I, \mathcal{B}_{n-2,n-2}^n)$. }\label{bbscan}
}
\end{figure}

The above two propositions are small variants of the arguments in \cite{BG}, so we leave them to the reader. Proposition
\ref{twobb} is obtained by a direct comparison, i.e. these two diffeomorphisms are induced by the same isotopy-extension construction. 

Proposition \ref{bbtriv} has an alternative way of being expressed.  Given the barbell diffeomorphism family, 
$$S^{i+j-n} \to \Diff(\mathcal{B}_{i,j}^n)$$
we can imagine this family fibering, i.e.
$$S^{i+j-n} \times D^{n-j-1} \to \Diff(\mathcal{B}_{i,j}^{j+1}).$$
Now consider the mid-ball in $\mathcal B_{i,j}^n$, this is a copy of $D^{n-1}$. The preimage the two cocores
in the mid-ball is given by the intersection of the map $S^{i+j-n} \times D^{n-1} \to \mathcal B_{i,j}^n$ with the
cocores, thus they will be two disjoint (framed) closed manifolds of dimension $(i-1)$ and $(j-1)$ respectively
in $S^{i+j-n} \times D^{n-1}$.  If we further use the fibering, we can imagine this as a $S^{i+j-n} \times D^{n-j-1}$-parametrized
family of $0$-manifolds and $(j-i)$-manifolds in the $\mathcal B_{i,j}^{j+1}$ mid-ball, which is a copy of $D^j$.  This 
family can be readily visualized.  The $0$-manifold family could be described as a parametrized family of null-cobordisms
of an embedded $S^0$, and the $(j-i)$-manifold family is similarly a parametrized null cobordism of $S^{j-i}$, i.e. 
this is a family of disjoint spheres: one a copy of $S^0$ and the other a copy of $S^{j-i}$ which on the boundary of
this $(i-1)$-dimensional family are spheres that bound disjoint discs -- which are used to construct the null cobordism. 
But since our family is $(i-1)$-dimensional this is exactly the right dimension that allows the family to link, which
is exactly what is going on.  Proposition \ref{bbtriv} is the homotopy-theoretic analogue of the linking number of this
parametrized family of high codimension links.  With a slight change of perspective we could perform this analysis in the
mid-ball of $\mathcal B_{i,j}^n$ using the $S^{i+j-n}$ parameter space.  This will be a family consisting generically
of two spheres: one $S^{n-j-1}$ and the other $S^{n-i-1}$ in $D^{n-1}$, thus a family parametrized by $S^{i+j-n}$
is precisely the correct dimension to allow for a linking number. 

While Proposition \ref{bbtriv} tells us that the inclusion $\Omega^{n-j} S^i \to \Diff(\mathcal{B}_{i,j}^n)$ is non-trivial
on the first non-trivial homotopy group, the inclusion is in fact a retract, i.e. non-trivial on all homotopy and homology
groups of $\Omega^{n-j} S^i$.  To show this we need to construct a map back. 

\begin{defn}
Observe there are maps
$$\Diff(\mathcal{B}_{i,j}^n) \to \Omega^{n-j} S^i, \hskip 1cm \Diff(\mathcal{B}_{i,j}^n) \to \Omega^{n-i} S^j$$
given by restricting to cocores and projecting to the cellular skeleton $S^i \vee S^j$, then forgetting the complementary sphere
wedge summand.  
\end{defn}

\begin{prop}\label{cocoreretract}
The inclusion of the fiber-preserving subspace $\Omega^{n-j} S^i \to \Diff(\mathcal{B}_{i,j}^n)$ is a retract, i.e. 
composition with the above map $\Diff(\mathcal{B}_{i,j}^n) \to \Omega^{n-j} S^i$ is homotopic to the identity. 
$$\xymatrix{\Omega^i S^i \ar[dr] \ar[rr]^{\simeq Id} & & \Omega^i S^i \\
               & \Diff(\mathcal{B}_{i,j}^n) \ar[ur] & }.$$

The composite with the other map
$$\xymatrix{\Omega^{n-j} S^i \ar[dr] \ar[rr]^{\simeq \Sigma^{(j-i)}} & & \Omega^{n-i} S^j \\
               & \Diff(\mathcal{B}_{i,j}^n) \ar[ur] & }$$
is homotopic, up to sign, to the iterated suspension map, i.e. $\Sigma^{j-i}$, i.e. we identify 
$\Omega^{j-i} S^i$ with the subspace of $\Omega^{n-i} S^j$ with the subspace of maps that preserve
$j-i$ suspension coordinates. 
\begin{proof}
The idea is to chase through the definition of our family, using the
fibration $\mathcal B_{i,j}^{j+1} \to \mathcal{B}_{i,j}^n \to D^{n-j-1}$. 
Thinking of the fiber as a once-punctured $S^i \times D^{j-i+1}$.  This gives us the inclusion 
$\Omega^{n-j} S^i \to \Diff(\mathcal{B}_{i,j})$ as fiber-preserving diffeomorphisms.   We consider the induced
diffeomorphisms of the $\mathcal{B}_{i,j}^{j+1}$ fibers.  The cocore complementary to the $S^i$ core sphere is
a copy of $D^{j+1-i}$, while the cocore complementary to the $S^j$ core sphere is a copy of the interval, $D^1$. 

The fact that the composite $\Omega^{n-j} S^i \to \Diff(\mathcal{B}_{i,j}) \to \Omega^{n-j} S^i$ is the identity map (after suitable
identifications) is derivable immediately from the definition, carefully keeping track of the suspension parameters. 

The composite $\Omega^{n-j} S^i \to \Diff(\mathcal{B}_{i,j}) \to \Omega^{n-i} S^j$ is depicted in Figure \ref{oisiojsj}.  The argument 
is essentially identical to the previous case, but our fibrewise cocores are copies of $D^{j+1-i}$, i.e. an interval with $j-i$
additional parameters.  These additional parameters supply the canonical null-homotopies of the embedded interval, which is
another way of stating that the map $\Omega^{n-j} S^i \to \Omega^{n-i} S^j$ is the suspension $\Sigma^{j-i}$. 

\begin{figure}[H]
{
\psfrag{bij}[tl][tl][1][0]{$\mathcal{B}_{i,j}^{j+1}$}
\psfrag{sidj1i}[tl][tl][0.7][0]{$S^i \times D^{j+1-i}$}
\psfrag{sjd1}[tl][tl][0.7][0]{$S^j \times D^1$}
$$\includegraphics[width=8cm]{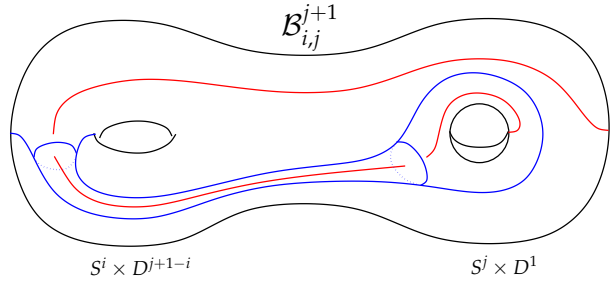}$$
}
\caption{\label{oisiojsj} The composite $\Omega^{n-j} S^i \to \Diff(\mathcal{B}_{i,j}) \to \Omega^{n-i} S^j$
in a fiber over a point in $D^{j-1}$. 
The image of $\{*\} \times D^{j+1-i}$ from the $S^i \times D^{j+1-i}$ summand is in blue, 
and the image of $\{*\} \times D^1$ from the $S^j \times D^1$ summand is in red.
}
\end{figure}

\end{proof}
\end{prop}


We  list one other elementary property of barbell diffeomorphisms. The idea is to consider the standard
inclusion $\Diff(\mathcal{B}_{i,j}^n) \to \Diff(S^i \times D^{n-i})$ and
$\Diff(\mathcal{B}_{i,j}^n) \to \Diff(S^j \times D^{n-j})$ obtained by attaching an $(i+1)$-handle or $(j+1)$-handle
respectively and extending via the identity map.  

\begin{prop}\label{nullhfillings}
The composites 
$$\Omega^{n-j} S^i \to \Diff(\mathcal{B}_{i,j}^n) \to \Diff(S^i \times D^{n-i})$$
and
$$\Omega^{n-j} S^i \to \Diff(\mathcal{B}_{i,j}^n) \to \Diff(S^j \times D^{n-j})$$
are canonically null-homotopic.
\begin{proof}
Recall the map $\Omega^{n-j} S^i \to \Diff(\mathcal{B}_{i,j}^n)$ was defined via a fibrewise isotopy-extension
process.  The first map in the statement of the proposition corresponds to forgetting the ball used to 
construct the isotopy extension, the second map corresponds to filling in the manifold in which the ball moves. 
In the first case, the diffeomorphism family is tautologically null due to the Palais homotopy long exact sequence. 
The second map is null as the input isotopy is itself null, i.e. the parametrizing family
$\Omega^{n-j} S^i$ factors through the inclusion $\Omega^{n-j} S^i \to \Omega^{n-j} D^{i+1}$. 
\end{proof}
\end{prop}

The rationale behind constructing the above null isotopies is that we can use them 
to construct certain null pseudo-isotopies, once we embed the barbell manifolds in larger manifolds.  This is
the content of Proposition \ref{nullPI}.



While the barbell diffeomorphisms themselves $\Omega^{n-j} S^i \to \Diff(\mathcal{B}_{i,j}^n)$ are not null in pseudo-isotopy, i.e.
they do not lift to maps $\Omega^{n-j} S^i \to \PDiff(\mathcal{B}_{i,j}^n)$, the implanted barbell diffeomorphisms are often null
in pseudoisotopy. The next proposition is a variation of Proposition \ref{nullhfillings}. 

\begin{prop}\label{nullPI}
Given an embedded barbell $\mathcal{B}_{i,j}^n \to N$ where $N$ is an $n$-manifold, the induced map 
$$\Omega^{n-j} S^i \to \Diff(N)$$
is null in pseudo-isotopy, provided one of the two core spheres is smoothly slice, i.e. is the boundary of
a smoothly-embedded $D^{i+1}$ or $D^{j+1}$ in $N \times I$. Precisely, there is a lift of the barbell diffeomorphism
family
$$\xymatrix{ & \PDiff(N) \ar[d] \\ \Omega^{n-j} S^i \ar[ur] \ar[r] & \Diff(N)}.$$
\begin{proof}
The group $\PDiff(N)$ is the collection of all diffeomorphisms of $N \times I$ which restrict to the identity on
$N \times \{0\}$ and $(\partial N) \times I$, often called the group of {\it pseudo-isotopy diffeomorphisms.} 
The idea is to consider the
manifold $S^i \times D^{n-i}$ (or $S^j \times D^{n-j}$) as the barbell manifold $\mathcal{B}_{i,j}^n$ union an
$i+1$ (or $j+1$)-handle respectively.  We embed $\mathcal{B}_{i,j}^n \times I$ into $N \times I$ using the map
$f(p,t) = (g(p), t/2)$ where $g : \mathcal{B}_{i,j}^n \to N$ is our barbell embedding.  We embed the $(i+1)$ or
$(j+1)$-handle in $N \times I$ so that its intersection with $N \times [0, \frac{1}{2}]$ exists in $U \times [0,\frac{1}{2}]$
where $U$ is a small neighbourhood of 
$g(\mathcal{B}_{i,j}^n)$ in $N$.  We can do this by ensuring the height function for the smooth slice disc has height $> \frac{1}{2}$
outside of a small neighbourhood of the slice sphere. This ensures the handle, in its interior, is disjoint from the image of $f$.
The image of $f$ union this handle is diffeomorphic to $S^i \times D^{n-i}$ or $S^j \times D^{n-j}$ respectively, thus our family
of diffeomorphisms $\Omega^{n-j} S^i \to \Diff(N)$ extends to a diffeomorphism of $N \times I$, using the null-isotopy
of Proposition \ref{nullhfillings} on the image of $f$ union the handle, which extends to $N \times I$ via the identity map. 
\end{proof}
\end{prop}

Proposition \ref{nullPI} was inspired by a conversation with David Gay, who has alternative descriptions
of such null pseudoisotopies. 

\begin{figure}[H]
{
\psfrag{NxI}[tl][tl][0.7][0]{$N \times I$}
\psfrag{BxI}[tl][tl][0.7][0]{\textcolor{blue}{$f(\mathcal{B}_{i,j}^n \times [0,1])$}}
$$\includegraphics[width=8cm]{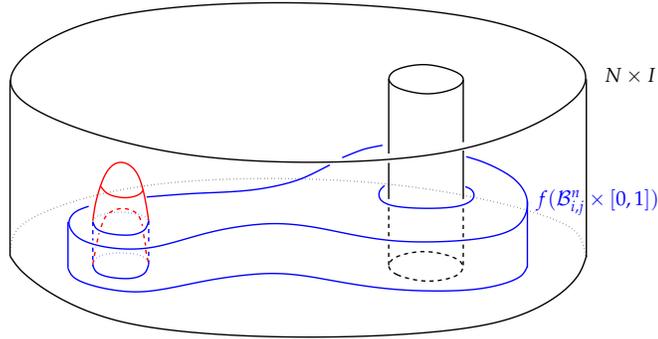}$$
\caption{\label{nullpifig} Null-pseudoisotopy via embedded null isotopy. Embedded handle in red.}
}
\end{figure}


We give a surgery description of the barbell diffeomorphisms in the case $i+j=n$.  
We start with the observation that one full
Dehn twist about about a curve in a punctured disc can be visualized by a technique of embedded surgeries. 

\begin{figure}[H]
{
\psfrag{km1}[tl][tl][0.7][0]{$k-1$}
$$\includegraphics[width=14cm]{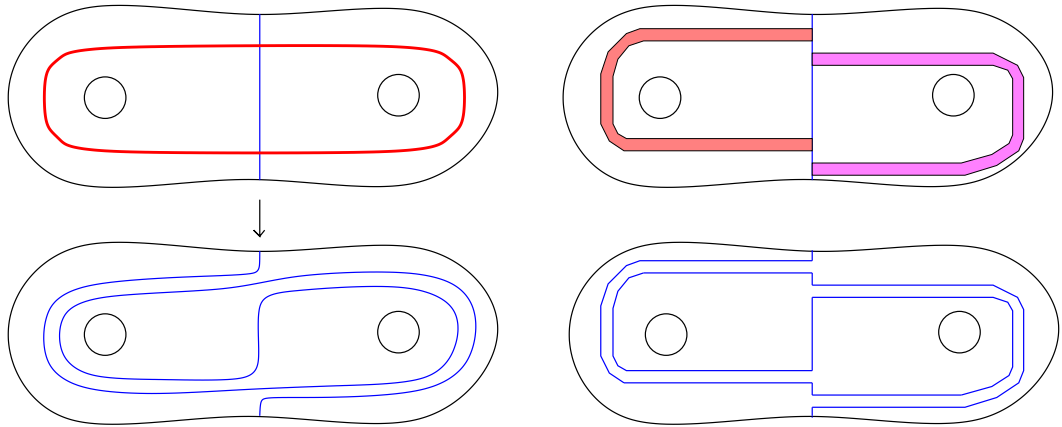}$$
\caption{\label{surjdehn}Surgery description of a Dehn twist}
}
In the upper-left figure we see a blue arc splitting the twice-punctured disc into
two annuli.  We perform a Dehn twist about the red circle, with the resulting embedded arc
appearing in the bottom-left.  In the top right we have two linking copies of $S^0$ embedded
in the blue arc, representing the attaching maps for two one-handles on the left (in orange) and
right (in magenta).  The result of the embedded surgery appears in the bottom-right. 
\end{figure}

The barbell diffeomorphism of $\mathcal{B}_{i,j}^n$ for $i+j=n$ has an analogous description.  One replaces the
blue arc in Figure \ref{surjdehn} by the mid-ball (diffeomorphic to $D^{n-1}$).  And one replaces the orange and
magenta $1$-handle attachments with $i$ and $j$-handle attachments respectively, with the $i$-handle being
the core of the $S^i \times D^j$ summand, and the $j$-handle attachment being the core of the $S^j \times D^i$
summand of $\mathcal{B}_{i,j}^n$.  The important issue is that the boundaries of the handle attachments 
are {\it linked} spheres in the mid-ball $S^{i-1} \sqcup S^{j-1} \to D^{n-1}$. 

\begin{prop} \label{handle_prop}
The action of the barbell diffeomorphism on the mid-ball of $\mathcal B_{i,j}^n$ when $n=i+j$ is isotopic
to replacing the mid-ball by its surgered embedding, where one does surgery on a trivially framed
link $S^{i-1} \sqcup S^{j-1} \subset D^{n-1}$ (the mid-ball) where the first sphere is attaching map for the
 core of the $S^i \times D^{n-i}$ summand, and the $S^{j-1}$ is the attaching sphere for the core of the
 $S^j \times D^{n-j}$ summand.  The link $S^{i-1} \sqcup S^{j-1} \subset D^{n-1}$ has unknotted components, 
 but the components have linking number $\pm 1$. 
\begin{proof}
To see this, consider $\mathcal B_{i,j}^n$ fibering over $D^{i-1}$ with fiber
$\mathcal B_{i,j}^{j+1}$.  Consider the action of the barbell diffeomorphism on the mid-balls in
the fibers. Specifically, consider the intersection of the image of these mid-balls with the cocores.
Generally these will consist of a disjoint union $S^0 \sqcup S^{j-i}$.  The $S^0$ comes from the $S^j$
cocore, while the $S^{j-i}$ comes from the $S^i$ cocore.  At the centre of the $D^{i-1}$ parameter space
the $S^{j-i}$ and $S^0$ sit on a common $D^{j-i+1}$ with one point of $S^0$ inside the $S^{j-i}$ and the 
other on the outside. As one moves the $D^{i-1}$ parameter the $S^0$ is pushed out of the subspace of the $S^{j-i}$, 
and as one approaches the boundary first the $S^{j-i}$ is coned-off, then the $S^0$ is coned-off.  This is exactly
the slicing perspective on the standard linked pair $S^{j-1} \sqcup S^{i-1} \subset D^{n-1}$, slicing over $D^{i-1}$. 
\end{proof}
\end{prop}

The advantage of this perspective is that it allows us to give a relatively elementary combinatorial description
of the barbell diffeomorphism, in terms of handle attachments.  

\begin{prop}\label{bbpi} The barbell diffeomorphism, as an element of $\pi_0 \Diff(\mathcal{B}_{i,j}^n)$ with $i+j=n$
is non-trivial in pseudo-isotopy.  We have two arguments.
The restriction to the mid-ball 
$$\pi_0 \Diff(\mathcal{B}_{i,j}^n) \to \pi_0 \Emb(D^{n-1}, \mathcal{B}_{i,j}^n)$$
is non-trivial in pseudo-isotopy, indeed if we let $\Map (D^{n-1}, \mathcal{B}_{i,j}^n)$ denote the space of maps of 
$D^{n-1}$ to $\mathcal{B}_{i,j}^n$ that restrict to the standard inclusion (the boundary connect-sum splitting disc) on the
boundary, then the map $\pi_0 \Diff(\mathcal{B}_{i,j}^n) \to \pi_0 \Map(D^{n-1}, \mathcal{B}_{i,j}^n)$ is homotopically non-trivial.
This latter space, up to a canonical homotopy-equivalence, is $\Omega^{n-1}(S^i \vee S^j)$.

The restriction to either cocore
$$\pi_0 \Diff(\mathcal{B}_{i,j}^n) \to \pi_0 \Emb(D^i, \mathcal{B}_{i,j}^n) \hskip 1cm \text{ or } \hskip 1cm 
  \pi_0 \Diff(\mathcal{B}_{i,j}^n) \to \pi_0 \Emb(D^j, \mathcal{B}_{i,j}^n)$$
is non-trivial in pseudo-isotopy.  Similarly, if we go one step further, the map 
$\pi_0 \Diff(\mathcal{B}_{i,j}^n) \to \pi_0 \Map(D^i, \mathcal{B}_{i,j}^n)$ is non-trivial.
\begin{proof}
The key observation is that the barbell diffeomorphism restricted to the mid-ball is obtained by surgery on a $2$-component link, 
$S^{i-1} \sqcup S^{j-1} \subset D^{n-1}$ corresponding to the core $S^i$ and $S^j$ respectively, i.e. the $S^{i-1}$ is the attaching sphere
for the $i$-handle, and the $S^{j-1}$ is the attaching sphere for the $j$-handle, when building $\mathcal{B}_{i,j}^n$ from the midball
by handle attachments.  Given an embedding $D^{n-1} \to \mathcal{B}_{i,j}^n$ it induces an element
of $\Omega^{n-1} \mathcal{B}_{i,j}^n \simeq \Omega^{n-1} (S^i \vee S^j)$ and the barbell diffeomorphism induces the
Whitehead product $[w_i, w_j]$ where $w_i : S^i \to S^i \vee S^j$ is the inclusion of $S^i$, and $w_j : S^j \to S^i \vee S^j$ 
is the inclusion of $S^j$. 

A second pseudo-isotopy obstruction follows from Proposition \ref{cocoreretract}.  Specifically, 
the embedding of the $i$-dimensional cocore may be projected to the $S^i$-core, giving an element of $\Omega^i S^i$. 
For the barbell diffeomorphism, this is a generating element of $\pi_0 \Omega^i S^i \equiv \pi_i S^i \simeq \Zed$.
\end{proof}
\end{prop}


\begin{prop}\label{bbtubing}
The action of the barbell diffeomorphism for $i+j=n$ on the cocores corresponds to tubing with the complementary
core sphere. The intersection of the image of the cocores with the mid-balls are Hopf-linked 
embedded copies of $S^{i-1} \sqcup S^{j-1}$, as in the Figure \ref{cocore-tubing}.
\begin{figure}[H]
{
\psfrag{si}[tl][tl][0.7][0]{$S^i$}
\psfrag{sj}[tl][tl][0.7][0]{$S^j$}
\psfrag{sim1}[tl][tl][0.7][0]{$S^{i-1}$}
\psfrag{sjm1}[tl][tl][0.7][0]{$S^{j-1}$}
\psfrag{di}[tl][tl][0.7][0]{$D^i$}
\psfrag{dj}[tl][tl][0.7][0]{$D^j$}
$$\includegraphics[width=14cm]{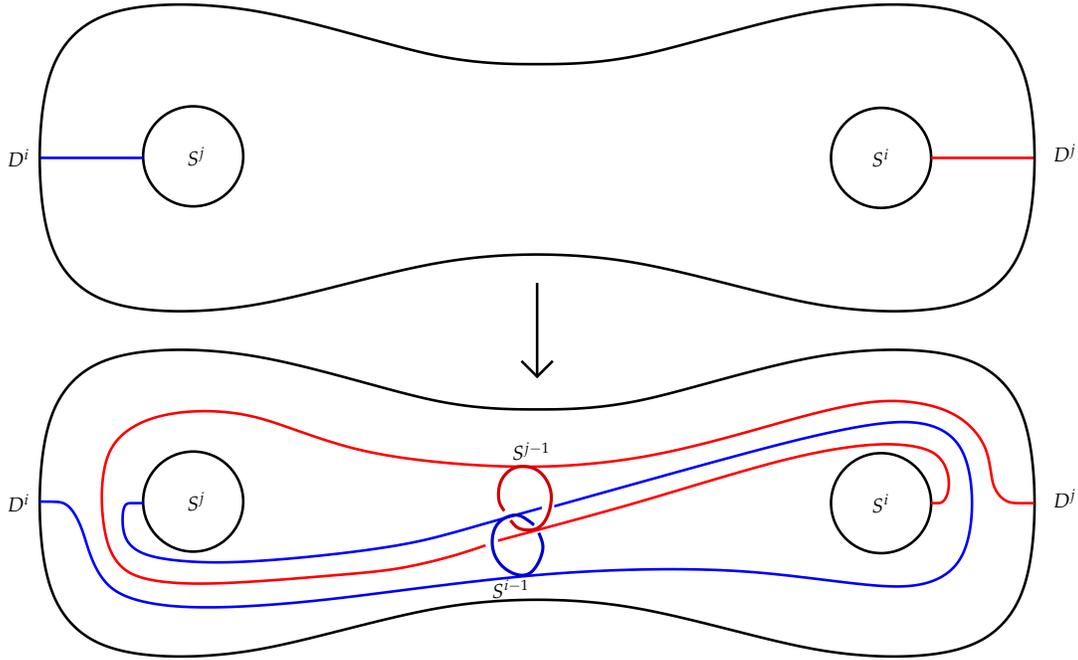}$$
\caption{Barbell applied to cocores having `linked' tubings.}\label{cocore-tubing}
}
\end{figure}
\begin{proof}
Consider $\mathcal B_{i,j}^n$ fibering over $D^{i-1}$ with fiber
$\mathcal B_{i,j}^{j+1}$.  The map $D^{i-1} \to \Diff( \mathcal B_{i,j}^{j+1} )$
corresponds to the diffeomorphisms induced by ambient isotopy as one slides the
`$j$ puncture' about the `$i$ puncture'.  The cocores in the fiber $\mathcal{B}_{i,j}^{j+1}$ 
correspond to embedded copies of $D^1$ (corresponding to the $j$-puncture) and $D^{j-i+1}$ (the $i$-puncture) respectively.
In the $i=1$ case, the cocore pair are linked as described, by an explicit performance of the 
isotopy-extension.  When $i>1$ the cocores are no longer linked in $\mathcal{B}_{i,j}^{j+1}$.  What we see is a fibering 
of the standard linked pair, fibered over $D^{i-1}$. 
\end{proof}
\end{prop}

Notice we have several equivalent ways to distinguish the barbell diffeomorphism from its inverse.  
Proposition \ref{bbpi} tells us that if we consider the intersection of the mid-ball with the image of the 
cocores, we get a standard linked pair.  If we orient the linked pair using normal bundles (i.e. the standard 
in oriented intersection theory) this would be a labeled and oriented 2-component link.  The linking number
is therefore a well-defined integer and these will be opposite for the barbell diffeomorphism and its inverse. 

\begin{figure}[H]
{
\psfrag{si}[tl][tl][0.7][0]{$S^i$}
\psfrag{sj}[tl][tl][0.7][0]{$S^j$}
\psfrag{sim1}[tl][tl][0.7][0]{$S^{i-1}$}
\psfrag{sjm1}[tl][tl][0.7][0]{$S^{j-1}$}
\psfrag{di}[tl][tl][0.7][0]{$D^i$}
\psfrag{dj}[tl][tl][0.7][0]{$D^j$}
$$\includegraphics[width=14cm]{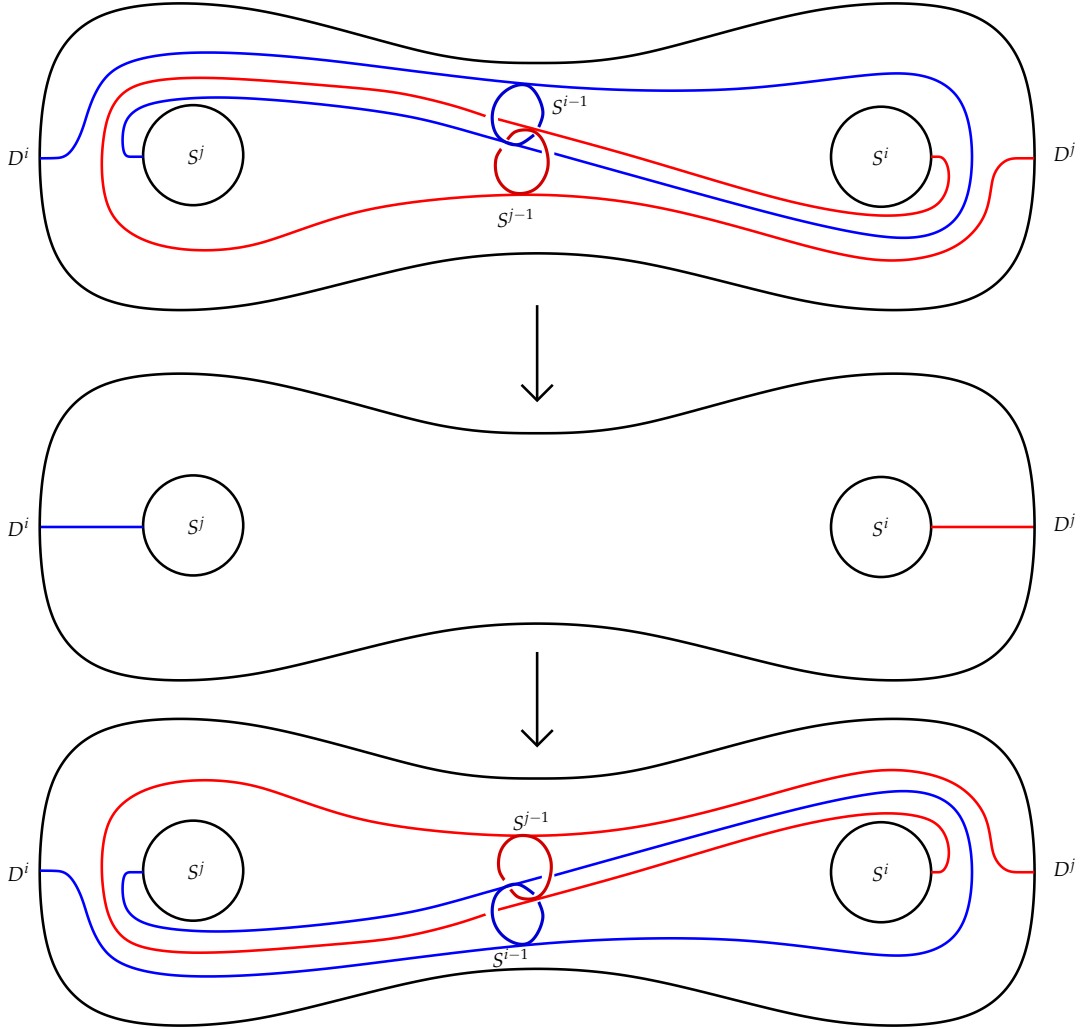}$$
\caption{Barbell image and pre-image of cocores.}\label{cocore-tubing2}
}
\end{figure}

There is an analogous result when $i+j>n$.  In this case, the cocores intersected with the mid-ball are too low-dimensional
to link, but in the {\it family} of maps $S^{i+j-n} \to \Omega^{n-j} S^i \to \Diff(\mathcal{B}_{i,j}^n)$ the family of cocores
intersected with the mid-ball analogously link. 

Barbell diffeomorphisms are closely related to the diffeomorphisms constructed by Watanabe \cite{W}.  See \cite{BG} for details. 


\section{$\pi_{n-4} \Diff(S^1 \times D^{n-1})$ and the $\delta_k$ diffeomorphisms. }\label{deltk}

In this section prove Theorem \ref{mainsmooththm} which is a new proof that the homotopy group 
$\pi_{n-4} \Diff(S^1 \times D^{n-1})$ is not finitely generated for $n \geq 4$.  On the large scale, this proof has several 
similarities to the one presented in \cite{BG} in that we 
compute the same $W_3$-invariant on the same implanted barbell diffeomorphisms $\delta_k$ to show they are linearly independent.   
The principal difference between the argument 
given here, and the one in \cite{BG} is that the computation of the $W_3$-invariant given here is directly from our definition of the
invariant $W_3$ and diffeomorphisms $\delta_k$.  In \cite{BG} 
we deduced relationships between the $W_3$-invariants of `nearby' implanted barbell diffeomorphisms, somewhat like bilinearity or a Skein relation.  
This relationship gave 
us a tool to reduce the computation of the $W_3$-invariant of any implanted barbell diffeomorphism
with linearly-embedded cuffs to that of $W_3(G(p,q))$.

\begin{figure}[H]
{
\psfrag{s1}[tl][tl][1][0]{$S^1$}
\psfrag{1}[tl][tl][0.7][0]{$1$}
\psfrag{2}[tl][tl][0.7][0]{$2$}
\psfrag{3}[tl][tl][0.7][0]{$3$}
\psfrag{km3}[tl][tl][0.7][0]{$k-3$}
\psfrag{km2}[tl][tl][0.7][0]{$k-2$}
\psfrag{km1}[tl][tl][0.7][0]{$k-1$}
$$\includegraphics[width=10cm]{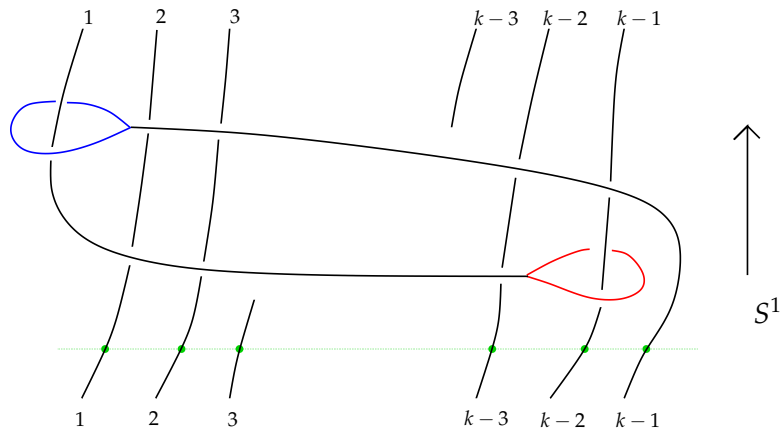}$$
\caption{$\delta_k$ barbell in $S^1 \times D^{n-1}$}\label{deltkfig}
}
\end{figure}

The elements $\delta_k \in \pi_{n-4} \Diff(S^1 \times D^{n-1})$ are the implanted barbells diffeomorphisms that come from
embeddings of the $\mathcal{B}_{n-2,n-2}^n$
barbells using `handcuff embeddings' as depicted in Figure \ref{deltkfig}.  The element of 
$\pi_{n-4} \Diff(\mathcal{B}_{n-2,n-2}^n)$ corresponds
to the image of the first non-trivial homotopy group ($\pi_{i+j-n}$), under the map 
$\Omega^{n-j} S^i \to \Diff(\mathcal{B}_{i,j}^n)$ when $i=j=n-2$ 
defined in Section \ref{bbdif}.

$$\Omega^2 S^{n-2} \to \Diff(\mathcal{B}_{n-2,n-2}^n).$$ 

\begin{figure}[H]
{
\psfrag{1}[tl][tl][0.7][0]{$1$}
\psfrag{2}[tl][tl][0.7][0]{$2$}
\psfrag{3}[tl][tl][0.7][0]{$3$}
\psfrag{km3}[tl][tl][0.7][0]{$k-3$}
\psfrag{km2}[tl][tl][0.7][0]{$k-2$}
\psfrag{km1}[tl][tl][0.7][0]{$k-1$}
$$\includegraphics[width=10cm]{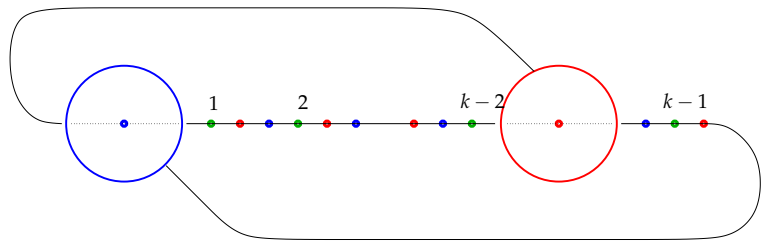}$$
\caption{Projection of $\delta_k$ barbell to $D^{n-1}$}\label{deltkproj}
}
\end{figure}

The inspiration for the $W_3$-invariant comes from Proposition \ref{bbtriv}, and it can also be seen in Figure \ref{cocore-tubing2}. 
Specifically, Proposition \ref{bbtriv} states that barbell diffeomorphisms are detectable by considering the mid-ball $B^{n-1}$ to
be fibered by intervals, giving a map $\Diff(\mathcal{B}_{n-2,n-2}^n) \to \Omega^{n-2} \Emb(I, \mathcal{B}_{n-2,n-2}^n)$.  In this formulation
we consider pairs of points $t_1 < t_2 \in I$ such that the embedding sends $t_1$ to the first cocore, and $t_2$ to the second, as a
signed intersection number for the family.  
Another way to state this is we are counting the linking number of the standard linking pair, of the pre-image of the cocores, intersected
with the mid-ball, i.e. it is a double-point formula for the linking number of the pair depicted in Figure \ref{cocore-tubing2}. 
Our preference is to state our invariant as a map of the form
$$\pi_{n-4} \Diff(\mathcal{B}_{n-2,n-2}^n) \to \pi_{n-4} \Omega^{n-2} \Emb(I, \mathcal{B}_{n-2,n-2}^n) \equiv \pi_{2n-6} \Emb(I, \mathcal{B}_{n-2,n-2}^n)\to \Zed$$
as this is an expression that we can generalize to $\Diff(S^1 \times D^{n-1})$. 
Since we understand the barbell diffeomorphism when restricted to the mid-ball, we can similarly `scan' through 
$\{1\} \times D^{n-1} \subset S^1 \times D^{n-1}$ giving a map 
$\Diff(S^1 \times D^{n-1}) \to \Omega^{n-2} \Emb(I, S^1 \times D^{n-1})$.  The $W_3$-invariant of an element
of $\pi_{n-4} \Diff(S^1 \times D^{n-1})$ takes
values in $\Rat \otimes \pi_{2n-6} \Emb(I, S^1 \times D^{n-1})$.  This homotopy-group is detectable at the $3^{rd}$-stage of
the Taylor tower, thus we consider the induced map of $3$-point configuration spaces to extract invariants of the map. 

\begin{figure}[H]
{
\psfrag{s1}[tl][tl][1][0]{$S^1$}
\psfrag{1}[tl][tl][0.7][0]{$1$}
\psfrag{2}[tl][tl][0.7][0]{$2$}
\psfrag{3}[tl][tl][0.7][0]{$3$}
\psfrag{4}[tl][tl][0.7][0]{$4$}
\psfrag{5}[tl][tl][0.7][0]{$5$}
\psfrag{6}[tl][tl][0.7][0]{$6$}
\psfrag{7}[tl][tl][0.7][0]{$7$}
\psfrag{8}[tl][tl][0.7][0]{$8$}
\psfrag{9}[tl][tl][0.7][0]{$9$}
\psfrag{10}[tl][tl][0.7][0]{$10$}
\psfrag{11}[tl][tl][0.7][0]{$11$}
\psfrag{12}[tl][tl][0.7][0]{$12$}
\psfrag{km4}[tl][tl][0.7][0]{$k-4$}
\psfrag{km3}[tl][tl][0.7][0]{$k-3$}
\psfrag{km2}[tl][tl][0.7][0]{$k-2$}
\psfrag{km1}[tl][tl][0.7][0]{$k-1$}
$$\includegraphics[width=10cm]{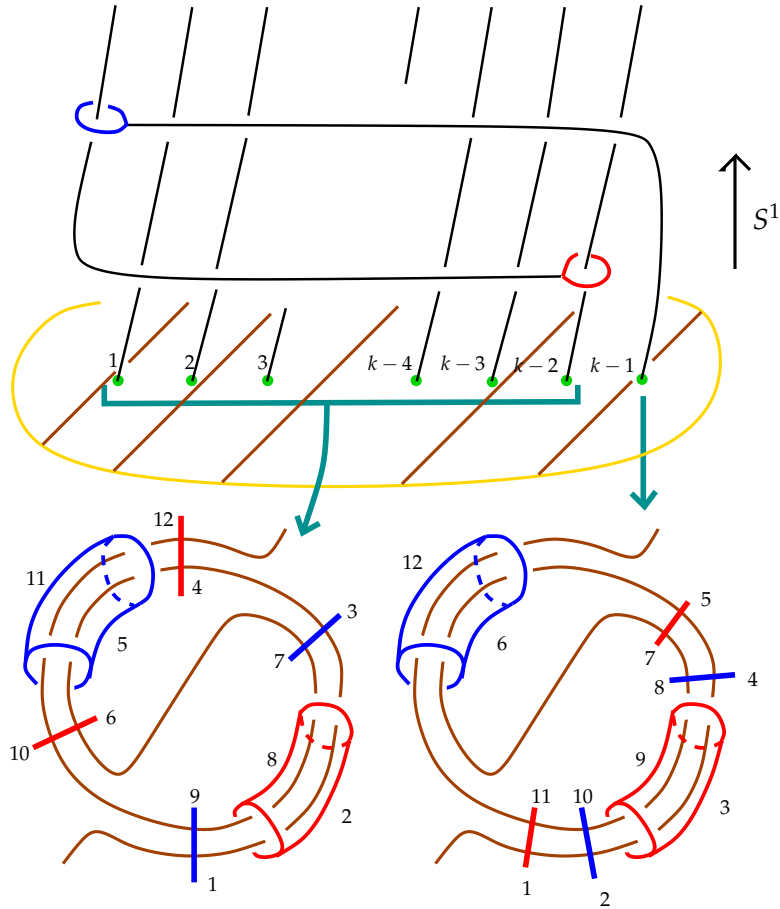}$$
\caption{Scanning $\delta_k$ barbell in $S^1 \times D^{n-1}$. Interval fibers
of mid-ball $\{1\} \times D^{n-1}$ in brown. }\label{twelvedots}
}
\end{figure}

Notice when scanning through $\delta_k$, if the brown interval fiber is disjoint from the bar, it is
unaffected by $\delta_k$.  But when it passes through the bar, imagine the bar's cross-section as 
$D^{n-1} \simeq D^{n-2} \times I$.  The
$D^{n-2}$ together with the $(n-4)$-parameter family of $\delta_k : S^{n-4} \to \Diff(S^1 \times D^{n-1})$ gives
us a $(2n-6)$-parameter family of embedded intervals in $S^1 \times D^{n-1}$, which are described in 
Proposition \ref{bbtriv} and Figure \ref{bbscan}. We modify Figure \ref{bbscan} as our barbell is embedded in
$S^1 \times D^{n-1}$ in handcuff fashion.  So as our brown interval passes through the bar, it curls as described in Figure \ref{twelvedots}. 
If we think of our $(2n-6)$-parameter family as $D^{n-3} \times D^{n-3}$, the first $D^{n-3}$ factor corresponds to 
the suspension parameter $S^{n-2} \equiv \Sigma S^{n-3}$ and controls the red cylinder being swung around the red
cuff. Similarly, the second copy of $D^{n-3}$ corresponds to the suspension parameter of the blue cuff $S^{n-2} \equiv \Sigma S^{n-3}$
and parametrizes the blue cylinder being swung around the blue cuff.  Lastly, the red and blue straight lines ($8$ in total)
depicted at the bottom of Figure \ref{twelvedots} indicate the points of the embedding that intersect the spanning disc for the
cuffs, i.e. the points on the embedding that have parameters with double-points. 

\begin{thm}\label{mainsmooththm}
$$W_3(\delta_k) = (k-1)\left(t_1^{-1} t_3^{1-k} + (-1)^{n} t_1^{1-k} t_3^{-1} - t_1^{2-k}t_3^1 + (-1)^{n-1} t_1 t_3^{2-k}\right) +$$
 $$ t_1t_3^{k-1} + (-1)^{n} t_1^{k-1} t_3 - t_1^{1-k} t_3^{2-k} + (-1)^{n-1} t_1^{2-k} t_3^{1-k}.$$
 The $W_3$-invariant takes values in the group $\Rat \otimes \pi_{2n-3} C_3'[S^1 \times D^{n-1}] / R$ and the
 elements $\{ W_3(\delta_k) : k \geq 4 \}$ are linearly-independent over $\Rat$.
\end{thm}

The remainder of this section is devoted to explaining the above: what precisely the
group $\Rat \otimes \pi_{2n-3} C_3'[S^1 \times D^{n-1}] / R$ is, how it can be considered a subgroup of 
$\Rat \otimes \pi_{2n-6} \Emb(I, S^1 \times D^{n-1})$, and how we compute $W_3(\delta_k)$ from 
Figure \ref{twelvedots} using only the double-point data.  
That said, the claimed formula above for $W_3(\delta_k)$ has some clear features in
common with Figure \ref{twelvedots}.  Notice that the bar crosses $\{1\} \times D^{n-1}$ at $(k-1)$ locations,
marked with green dots. There similarly a term with coefficient $k-1$ in the $W_3(\delta_k)$ formula.    Roughly speaking, the remaining 
term in the $W_3(\delta_k)$ computation is a correction term, since there is a different combinatorial pattern in the 
double-point data for the family as it crosses through the the green dot labelled $k-1$.

\begin{defn}If $M$ is a manifold, the {\bf configuration space} of $k$ points in $M$ is the
space
$$C_k(M) = \{ (p_1,\cdots,p_k) : p_i \neq p_j \ \forall i \neq j\}.$$
The Fulton-Macpherson compactification of $C_k(M)$ is denoted $C_k[M]$. This is obtained by 
taking the closure of $C_k(M)$ under the product map 
$$C_k(M) \to M^k \times (S^n)^{k \choose 2} \times [0,\infty]^{k \choose 3}$$
where the inclusion $C_k(M) \to M^k$ is set-theoretic inclusion $C_k(M) \subset M^k$.  The 
maps $C_k(M) \to S^n$ come from taking unit displacement vectors between pairs of points, $\frac{p_i-p_j}{|p_i-p_j|}$
i.e. we
assume $M \subset \Real^{n+1}$. Lastly, the maps  $C_k(M) \to [0, \infty]$ come from 
the relative ratio map $\frac{|p_i-p_j|}{|p_i-p_l|}$ where $\{i,j,l\} \subset \{1,2,\cdots,k\}$.
\end{defn}

Provided $M$ is compact $C_k[M]$ is a compact manifold with corners.  Moreover, 
the construction is natural with respect to embeddings and the inclusion $C_k(M) \to C_k[M]$ is a 
homotopy-equivalence.

In the context where we are considering embedding spaces $\Emb(I, M)$, the notation $C_k'[M]$ indicates a small
variation of $C_k[M]$ where.  Specifically, we take the pull-back of the unit tangent bundle under the map
$$\xymatrix{ & UTM^{k+2} \ar[d] \\ C_{k+2}[M] \ar[r] & M^{k+2}}$$
and then restrict to the subspace where $p_1$ (and its vector) agree with the initial-point 
of the embeddings of $\Emb(I, M)$, and $p_{k+2}$ (and its vector) agree with the terminal point of the embeddings 
defining $\Emb(I, M)$.  See \cite{Sin2} for details.   In the case of the interval, we restrict to the 
subspace of $C_k'[I]$ such that $t_1 \leq t_2 \leq \cdots \leq t_k$, i.e. we choose a standard connected-component. 
In this case, $C_k'[I]$ is known to be the $k$-th Stasheff Polytope, or associahedron.

Denote the generators of $\pi_1 C_k[S^1 \times D^{n-1}] \simeq \Zed^k$ by $\{t_i : i = 1,2,\cdots,n\}$. 
The class $w_{ij} \in \pi_{n-1} C_k[S^1 \times D^{n-1}]$ has all $k$ points stationary, 
with the exception of point $j$ that orbits around point $i$. 

\vskip 0.4cm

$\pi_{n-1} C_k[S^1 \times D^{n-1}]$ is generated by the set $\{t_l^q.w_{ij} \ \forall i,j,l,q\}$, with the relations

\begin{itemize}
\item $w_{ii} = 0 \ \forall i$
\item $w_{ij} = (-1)^n w_{ji} \ \forall i \neq j$. 
\item $t_l.w_{ij} = w_{ij}$ provided $l \notin \{i,j\}$.
\item $t_j.w_{ij} = t_i^{-1}.w_{ij} \ \forall i,j$.   
\end{itemize}
 
The way one proves the above is to observe the forgetful map $C_k(S^1 \times D^{n-1}) \to C_{k-1}(S^1 \times D^{n-1})$
is a locally-trivial fiber bundle.  Moreover, the map has a section, so the homotopy-groups of $C_k(S^1 \times D^{n-1})$
are isomorphic to the product of the homotopy groups of the fibers (iteratively), which are (individually) wedges of $S^1$ with 
various copies of $S^{n-1}$. The $S^1$ factors contribute the $t_i$ generators in $\pi_1$, while the sphere factors 
contribute the $w_{ij}$ generators.  By the Hilton-Milnor theorem the higher rational homotopy groups are generated
by Whitehead Products.  The Whitehead Product is a bilinear mapping $[\cdot,\cdot] : \pi_i X \times \pi_j X \to \pi_{i+j-1} X$
satisfying the 
$$(-1)^{pr}[[f,g],h] + (-1)^{pq}[[g,h],f] + (-1)^{rq}[[h,f],g] = 0,$$
$$\text{where } f \in \pi_p X, g \in \pi_q X, h \in \pi_r X \text{ with } p,q,r \geq 2.$$
Due to the form of the above relation it is sometimes called a `graded Jacobi identity' in analogy with the Lie Bracket. 

There are two elementary relations satisfied by the $w_{ij}$ classes via the Whitehead product:

\begin{itemize}
\item $[w_{ij},w_{lm}] = 0$ when $\{i,j\} \cap \{l,m\} = \emptyset$. 
\item $[w_{ij}+w_{il}, w_{lj}] = 0$ for all $i,j,l$. 
\end{itemize}

The latter relation should be viewed a generalized `orbital system' map $S^n \times S^n \to C_3(D^n)$ where there is an 
earth-moon-sun orbital triple.  For this 
interpretation one views the Whitehead Bracket as the obstruction to extending a wedge of maps, i.e. $S^i \vee S^j \to X$ to the
product $S^i \times S^j \to X$.  For an orbital triple such a map exists, thus the corresponding Whitehead bracket is zero. 

The latter relation above can be rewritten as $[w_{ij}, w_{jk}] - [w_{jk}, w_{ki}] = 0$,
giving the equality of the three cyclic permutations,
$$[w_{ij}, w_{jk}] = [w_{jk}, w_{ki}] = [w_{ki}, w_{ij}].$$

\begin{proposition}\label{htpyck} The rational homotopy-groups of $C_k[S^1 \times D^{n-1}]$ 
are generated by the Whitehead products of 
the elements $t_l^m.w_{ij}$. These satisfy the relations
\begin{itemize}
\item $[w_{ij}, w_{lm}]=0$ if $\{i,j\} \cap \{l,m\} = \emptyset$,
\item $[w_{ij}, w_{jl}] = [w_{jl}, w_{li}] = [w_{li}, w_{ij}]$,
\item $t_l.[f,g] = [t_l.f,t_l.g]$.
\end{itemize}
\end{proposition}

A relatively constructive way to verify much of the above is via intersection theory.  
Fix a unit  direction $\zeta \in \partial D^{n-1}$. Define $t^i Co_1^2$ to consist of pairs of points 
$(p_1,p_2) \in C_2(\Real^1 \times D^{n-1})$ such that the displacement vector $t_2^i.p_2 - p_1$ is a positive multiple of $\zeta$. 
We call $t^i Co_1^2$ a {\it cohorizontal manifold}.  Given an element of $\pi_{n-1} C_2(S^1 \times D^{n-1})$, we lift the
map to the universal cover $\tilde C_2(S^1 \times D^{n-1}) \subset C_2(\Real \times D^{n-1})$ and take its intersection
with the $t^i Co_1^2$ submanifold is a well-defined framed $0$-dimensional manifold as a cobordism class, thus an integer.  This 
invariant detects the class $t_1^i w_{12}$.  

To similarly detect homotopy classes in $\pi_{2n-3} C_3(S^1 \times D^{n-1})$ we have the collinear classes.  This will consist of three
points sitting on a `straight line' in $S^1 \times D^{n-1}$.  Roughly speaking by `straight line' we are referring to geodesics in the
standard Euclidean metric on $S^1 \times D^{n-1}$.   To be more precise, the manifold $Col^1_{\alpha, \beta}$ is the collection 
of points of the form $(p_1,p_2,p_3) \in C_3(\Real \times D^{n-1})$ such that $(p_2,t_1^\alpha p_1, t_3^\beta p_3)$ sit on a straight line
in $\Real \times D^{n-1}$ in the listed order. The manifold $Col^3_{\alpha, \beta}$ is similarly 
defined by the requirement $(t_1^\alpha p_1, t_3^\beta p_3, p_2)$ sit on a straight line
in $\Real \times D^{n-1}$ in the listed order.  The universal cover of $C_3(S^1 \times D^{n-1})$ is naturally an open subspace of
$(\Real^1 \times D^{n-1})^3$, thus we can consider $Col^1_{\alpha, \beta}$ and $Col^3_{\alpha, \beta}$ naturally as subspaces of the universal
cover of $C_3(S^1 \times D^{n-1})$. The manifolds $Col^1_{\alpha,\beta}$ and $Col^3_{\alpha,\beta}$
are disjoint and closed in the universal cover of $C_3(S^1 \times D^{n-1})$.   
Given a map $S^{2n-3} \to C_3(S^1 \times D^{n-1})$, we take its lift to the universal cover
$S^{2n-3} \to \tilde C_3(S^1 \times D^{n-1})$ and take the pre-image of
the pair $(Col^1_{\alpha,\beta}, Col^3_{\alpha,\beta})$.  Generically, this gives us a disjoint pair of compact oriented 
manifolds of dimension $(n-2)$ in $S^{2n-3}$, thus they have a well-defined linking number.  This linking number detects 
the coefficient of $t_1^\alpha t_3^\beta [w_{12},w_{23}]$. 

Given that linking numbers of pre-images of the pair $(Col^1_{\alpha,\beta}, Col^3_{\alpha,\beta})$ can be difficult to visualize
and compute for a lift of an arbitrary map $S^{2n-3} \to C_3(S^1 \times D^{n-1})$, we describe an isotopy of the
pair $(Col^1_{\alpha,\beta}, Col^3_{\alpha,\beta})$ that converts the computation into something that is often more manageable. 
For $\epsilon \in \Real$ consider the diffeomorphism of $\Real^n$ given by 
$$P_\epsilon(x_1, x_2, \cdots, x_n) = \left(x_1, x_2, \cdots, x_{n-1}, x_n + \epsilon \sum_{i=1}^{n-1} x_i^2\right).$$
This diffeomorphism has the feature that it converts the $x_n=c$ hyperplanes into paraboloids, when $\epsilon \neq 0$, similarly 
it turns lines in the $x_n=c$ plane into parabolas, but on a line parallel to the $x_n$-axis the diffeomorphism acts by translation. 
Moreover, $P_{\epsilon_1} \circ P_{\epsilon_2} = P_{\epsilon_1 + \epsilon_2}$ and $P_\epsilon^{-1} = P_{-\epsilon}$.  
If we consider the $Col^i_{\alpha,\beta}$ manifolds to be submanifolds naturally defined in $C_3(\Real^n)$ (i.e. before we pull them 
back to $C_3(\Real \times D^{n-1})$), we can pull them back via the diffeomorphism $P_\epsilon$, and these will be manifolds of 
coparabolic triples (plus triples on the lines parallel to the $x_n$-axis). As $P_\epsilon$ is an orientation-preserving diffeomorphism, 
these manifolds when pulled-back to $\tilde C_3(S^1 \times D^{n-1})$ also detect the $t_1^\alpha t_3^\beta [w_{12}, w_{13}]$ classes. 
We denote the pull-backs of the collinear manifolds the coparabolic manifolds, i.e. 
$$Cop^i_{\alpha,\beta,\epsilon} = P_\epsilon^*(Col^i_{\alpha,\beta}).$$

Given a map $S^{2n-3} \to C_3(S^1 \times D^{n-1})$, we lift to the universal cover and take the pre-images of $Cop^i_{\alpha,\beta,\epsilon}$ for $i=1,3$.  Given that these are disjoint closed, oriented manifolds in the codomain, their pre-images are disjoint compact, oriented manifolds
in $S^{2n-3}$.  Generically we can assume these maps contain no cohorizontal triples, since the cohorizontal triple condition is
of codimension $2n-2 > 2n-3$.   Thus for if $\epsilon$ is large, and our triples of points are not approximating infinitesimal triples 
along the embedding, we can assume that the critical point of the parabola occurs outside the $\Real \times D^{n-1}$, and this 
critical point separates two of the three points in the parabolic triple.  Thus in the limit, the linking of the pre-image of the pair 
$(Cop^1_{\alpha,\beta,\epsilon}, Cop^3_{\alpha,\beta,\epsilon})$ is computable as the linking numbers of the pre-image of the pair
$(t^\alpha Co_2^1 - t^{\alpha-\beta} Co_3^1, t^{\beta-\alpha} Co_1^3 - t^\beta Co_2^3)$, as steep segments of parabolas approximate 
cohorizontal lines.  The reason for the signs, such as the minus sign in front of the $t^{\alpha-\beta} Co_3^1$ term is that when 
computing the signed intersection number of a parabolic triple, all the signs for the $Co_2^1$ and $Co_3^1$ are the same, with the
exception for the reversal of direction of the parabola.  

\begin{figure}[H]
{
\psfrag{zeta}[tl][tl][0.7][0]{$\zeta$}
\psfrag{2}[tl][tl][0.7][0]{$2$}
\psfrag{3}[tl][tl][0.7][0]{$3$}
\psfrag{km3}[tl][tl][0.7][0]{$k-3$}
\psfrag{km2}[tl][tl][0.7][0]{$k-2$}
\psfrag{km1}[tl][tl][0.7][0]{$k-1$}
$$\includegraphics[width=14cm]{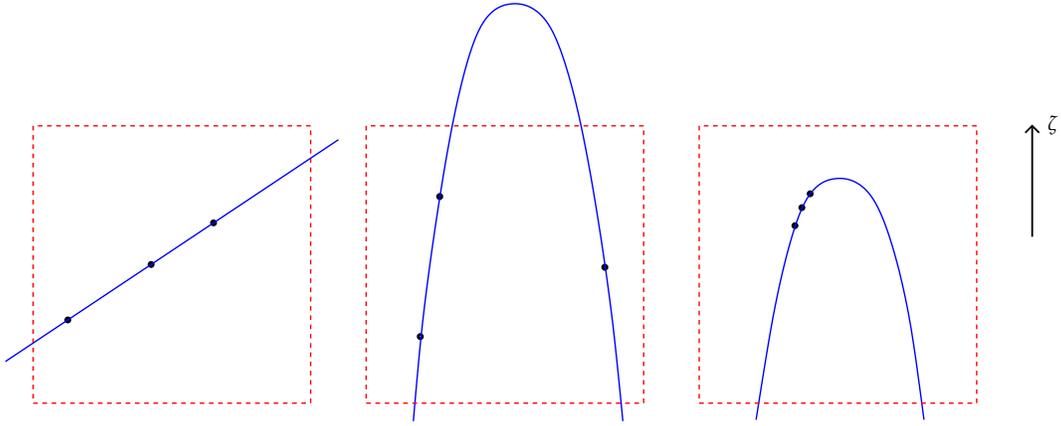}$$
\caption{Parabolic triples, $\epsilon=0$ left.  Large $\epsilon$ middle and right.}\label{parabdiag}
}
\end{figure}

\begin{lem}\label{polyakstyle}
Given a smooth map $f : S^{2n-3} \to C_3(S^1 \times D^{n-1})$, generically we can assume it has no cohorizontal triples.  Moreover, 
provided the map does not have any cohorizontal triples, consider the lift to the universal cover
$\tilde f : S^{2n-3} \to \tilde C_3(S^1 \times D^{n-1}).$ The linking numbers in $S^{2n-3}$ of the 
pre-images of the disjoint pair of manifolds in $\tilde C_3(S^1 \times D^{n-1})$
$$(Col^1_{\alpha,\beta}, Col^3_{\alpha,\beta})$$
agrees with the linking number of
the pre-image of the pair 
$$(t^\alpha Co_2^1 - t^{\alpha-\beta} Co_3^1, t^{\beta-\alpha} Co_1^3 - t^\beta Co_2^3)$$
of linear combinations of manifolds, for all $\alpha, \beta \in \Zed$.
\end{lem}

Lemma \ref{polyakstyle} is a variant of an argument the first author learned from Misha Polyak verbally in 2005, 
who later described his argument in more detail in his presentation \cite{PV}.
While Polyak's argument occurs at the $4^{th}$-stage of the Taylor tower, this variant works at the $3^{rd}$.  

As in \cite{BG}, given a map $S^{n-4} \to \Diff(S^1 \times D^{n-1})$ we compose with the scanning map
$\Diff(S^1 \times D^{n-1}) \to \Omega^{n-2} \Emb(I, S^1 \times D^{n-1})$.  We further compose with the 
$3^{rd}$ stage of the Taylor tower, which we think of as the evaluation map 
$\Emb(I, S^1 \times D^{n-1}) \times C_3'[I] \to C_3'[S^1 \times D^{n-1}]$, or its adjoint
$$ev_3 : \Emb(I, S^1 \times D^{n-1}) \to \Map( C_3'[I], C_3'[S^1 \times D^{n-1}])$$
where in this mapping space we demand that maps are stratum-preserving and aligned, meaning that when points 
collide in the domain, the corresponding points collide in the codomain, moreover, their associated tangent vectors
agree.  Putting these ingredients together we have an induced map
$$S^{n-4} \wedge S^{n-2} \equiv S^{2n-6} \to \Map( C_3'[I], C_3'[S^1 \times D^{n-1}]).$$

One could have some concerns here that the lifts of these maps to the universal cover of $C'_3[S^1 \times D^{n-1}]$
have transversality issues on various boundary strata of $C_3[I]$.  Specifically, if our family in $\Emb(I, S^1 \times D^{n-1})$
has any velocity vectors parallel to the `vertical' direction (i.e. the $S^1$-direction of $S^1 \times D^{n-1}$) then one
has collinear triples for infinitely many $\alpha, \beta$ of the $Col^i_{\alpha, \beta}$ variety.  For example, on the
$t_1=t_2=t_3$ stratum, but also there are intersections of distinct codimensions on the $t_1=t_2$ stratum, where the
co-dimension depends on the choice of $\alpha, \beta$.  There are several ways to avoid these transversality problems. 
For example, if the family has no vertical tangent vectors, this issue does not arise.  That said, vertical tangent vectors
are not known to be avoidable (that said, we do know the tangent vector field along the embeddings to be canonically null-homotopic).  
The underlying geometric problem is that the collinear manifolds contain all the vertical lines $\Real \times \{p\}$ for $p \in D^{n-1}$.  
Thus we can change our model to avoid this problem.  Using the coparabolic manifolds $Cop^i_{\alpha,\beta, \epsilon}$ 
with $\epsilon \neq 0$ suffices.
Our families $S^{2n-6} \times C_3'[I] \to C_3'[S^1 \times D^{n-1}]$ are generically transverse to these coparabolic manifolds, pulling them
back to oriented co-dimension $n-1$ submanifolds, along all strata. 

The associated map on the $2^{nd}$ stage $S^{2n-6} \to \Map( C_2'[I], C_2'[S^1 \times D^{n-1}])$ is
torsion (see \cite{BG} for details), so we can attach a null-homotopy to an appropriate multiple of the
$3^{rd}$ stage, giving a homotopy-class of map $S^{2n-3} \to C_3'[S^1 \times D^{n-1}]$.   Depending on 
which null-homotopy we attach, we can get a different homotopy-class of map
$S^{2n-3} \to C_3'[S^1 \times D^{n-1}]$. 
This is the subject of items (1)--(4) below. 

As we have seen $\pi_{2n-3} C_2'[S^1 \times D^{n-1}]$ is isomorphic to 
$\pi_{2n-3}(S^1 \vee S^{n-1}) \oplus \bigoplus_2 \pi_{2n-3} S^{n-1}$. 
Modulo torsion, the generators of $\pi_{2n-3}(S^1 \vee S^{n-1})$ are the Whitehead products of elements 
$t^k w_{12}$ for $k \in \Zed$.  
This gives us the result that $\pi_{2n-3} C_2'[S^1 \times D^{n-1}]$, mod torsion, is isomorphic to 
$\Zed[t_1^{\pm 1}, t_2^{\pm 1}] / \langle t_1t_2-1 = 0 \rangle$ as a module over the group-ring 
of the fundamental group. 
The generator of $\pi_{2n-3} C_2'[S^1 \times D^{n-1}]$ corresponding to a monomial $t_1^\alpha t_2^\beta$ is
$t_1^\alpha t_2^\beta w_{12}$. By attaching a homotopy-class of maps $S^{2n-6} \times I \times C_2[I] \to C_2'[S^1 \times D^{n-1}]$
to a closed-off $S^{2n-6} \times C_3[I] \to C_3'[S^1 \times D^{n-1}]$ we change the homotopy class by adding:

\begin{enumerate}
\item $[t_2^\alpha w_{23}, t_2^\beta w_{23}]$. This comes from the $t_1=0$ face. 
Thus the generator $t_1^\alpha w_{12}$ is mapped to $t_2^\alpha w_{23}$, and a Whitehead
bracket $[t_1^\alpha w_{12}, t_1^\beta w_{12}]$ is mapped to $[t_2^\alpha w_{23}, t_2^\beta w_{23}]$.

\item $[t_1^\alpha w_{12}, t_1^\beta w_{12}]$ to 
$[t_1^\alpha w_{13}+ t_2^\alpha w_{23} + a_1 w_{21}, t_1^\beta w_{13}+ t_2^\beta w_{23} + a_1 w_{21}]$. 
This comes from the $t_1=t_2$ face map, i.e. the inclusion $C_2'[S^1 \times D^{n-1}] \to C_3'[S^1 \times D^{n-1}]$
that doubles the first point, i.e. $(p_1, p_2) \longmapsto (p_1, \epsilon^+ p_1, p_2)$, where the 
perturbation $\epsilon^+ p_1$ is in the direction of the velocity vector. The integer $a_1$ is the degree of this
velocity vector map.   
This map sends $w_{12}$ to $w_{13}+w_{23} + a_1 w_{21}$,  $t_1$ to $t_1t_2$ and $t_2$ to $t_2$.  The 2nd stage
of the Taylor tower induces a null-homotopy of the velocity vector map, so we can assume $a_1=0$, but it is of interest
that the following computation gives the same answer for $a_1 \neq 0$. 
Thus it sends
$[t_1^\alpha w_{12}, t_1^\beta w_{12}]$ to $[t_1^\alpha w_{13}+ t_2^\alpha w_{23} + a_1 w_{21}, t_1^\beta w_{13}+ t_2^\beta w_{23} + a_1w_{21}]$.  Expanding this bracket using bilinearity we get
$$ = \left( -t_1^{\alpha-\beta} t_3^{-\beta} + (-1)^{n} t_1^{\beta-\alpha} t_3^{-\alpha} \right)[w_{12},w_{23}] + [t_1^\alpha w_{13}, t_1^\beta w_{13}] + $$
$$ a_1 \left( (-1)^{n-1} t_3^{-\beta} + (-1)^{n} t_3^{-\beta} + t_1^{-\alpha} - t_1^{-\alpha} \right)[w_{12},w_{23}]$$
where the latter row comes from collecting the terms involving $a_1$, and clearly these terms sum to zero. 

\item $[t_1^\alpha w_{12}+t_1^\alpha w_{13}+a_2 w_{23}, t_1^\beta w_{12}+ t_1^\beta w_{13}+a_2w_{23}]$. 
This is for the $t_2=t_3$ facet.  This corresponds to the map $C_2'[S^1 \times D^{n-1}] \to C_3'[S^1 \times D^{n-1}]$ that
doubles the second point, i.e. $(p_1,p_2) \longmapsto (p_1, p_2, \epsilon^+ p_2)$.  This map sends
$w_{12}$ to $w_{12} + w_{13} + a_2w_{23}$, $t_1$ to $t_1$ and $t_2$ to $t_2t_3$.  Thus
$[t_1^\alpha w_{12}, t_1^\beta w_{12}] \longmapsto [t_1^\alpha w_{12}+t_1^\alpha w_{13} + a_2w_{23}, t_1^\beta w_{12}+ t_1^\beta w_{13}+a_2w_{23}]$.  Like the previous case, this simplifies to
$$ = \left( -t_1^\alpha t_3^{\alpha-\beta} + (-1)^{n} t_1^\beta t_3^{\beta-\alpha} \right)[w_{12},w_{23}] + [t_1^\alpha w_{13}, t_1^{\beta}w_{13}] + $$
$$ a_2 \left( t_1^\beta - t_1^\beta + (-1)^{n-1} t_1^\alpha + (-1)^{n} t_1^\alpha \right)[w_{12},w_{23}].$$
Again, the terms with $a_2$ cancel.  

\item $[t_1^\alpha w_{12}, t_1^\beta w_{12}]$. This is for the $t_3=1$ facet.  
This corresponds to the inclusion $C_2'[S^1\times D^{n-1}] \to C_3'[S^1\times D^{n-1}]$ that
maps $(p_1,p_2)$ to $(p_1,p_2,(1,0))$, thus it sends $w_{12} \longmapsto w_{12}$ and $t_1 \longmapsto t_1$, 
$t_2 \longmapsto t_2$, thus it acts trivially on $[t_1^\alpha w_{12}, t_1^\beta w_{12}]$. 
\end{enumerate}

Thus our invariant via closure $\frac{1}{m}\overline{ev_3}(mf)$ of $\pi_{2n-6} \Emb(I, S^1 \times D^{n-1})$ takes values in
$$\Rat \otimes \pi_{2n-3} C_3'[S^1 \times D^{n-1}] / R$$
where $R$ is the subgroup generated by the above four inclusions.  Notice (1) kills the summand corresponding to
the $w_{23}$ brackets, and (4) kills the summands corresponding to the $w_{12}$ brackets. Using relation (1) and (4)
we can simplify (2) and (3) into relations between $w_{13}$ brackets and brackets of the form $[w_{12}, w_{23}]$, giving us
the Proposition \ref{tor-clo-arg2}. 

\begin{proposition}({\bf Closure Argument}) \label{tor-clo-arg2}
Given an element of 
$[f] \in \pi_{2n-6} \Emb(I, S^1 \times D^{n-1})$ such that $ev_2(f) : S^{2n-6} \to T_2 Emb(I, S^1 \times D^{n-1})$ is null,
we form the closure of the evaluation map $ev_3(f) : S^{2n-6} \to T_3 \Emb(I, S^1 \times D^{n-1})$ which is
a based map of the form
$$\overline{ev_3}(f) : S^{2n-3} \to C_3'[S^1 \times D^{n-1}].$$
The homotopy-class of this map, as a function of the homotopy-class $[f]$ is well-defined modulo a subgroup we call $R$.  
$R$ is generated by the torsion subgroup of
$\pi_{2n-3} C_3'[S^1 \times D^{n-1}]$ together with the elements 
$$\left(t_1^{\alpha-\beta}t_3^{-\beta} - t_1^\alpha t_3^{\alpha-\beta} + (-1)^{n} \left( t_1^\beta t_3^{\beta-\alpha} - t_1^{\beta-\alpha} t_3^{-\alpha} \right)\right) [w_{12},w_{23}] \ \forall \alpha,\beta \in \Zed, $$
$$ [t_2^\alpha w_{23}, t_2^\beta w_{23}] \ \forall \alpha, \beta, $$
$$ [t_1^\alpha w_{12}, t_1^\beta w_{12}] \ \forall \alpha, \beta, $$
$$ [t_1^\alpha w_{13}, t_1^\beta w_{13}] + \left(t_1^{\alpha-\beta}t_3^{-\beta} + (-1)^{n-1} t_1^{\beta-\alpha}t_3^{-\alpha}\right)[w_{12},w_{23}] \ 
\forall \alpha, \beta.$$

Since $\pi_{2n-6} T_2 \Emb(I, S^1 \times D^{n-1})$ is torsion, there is a homomorphism, called the {\bf closure operator}
$$\pi_{2n-6} \Emb(I, S^1 \times D^{n-1}) \to \Rat[t_1^{\pm 1}, t_3^{\pm 1}] / 
\langle 
t_1^{\alpha-\beta}t_3^{-\beta} - t_1^\alpha t_3^{\alpha-\beta} = (-1)^{n-1} \left( t_1^\beta t_3^{\beta-\alpha} - t_1^{\beta-\alpha} t_3^{-\alpha} \right)\ \forall \alpha,\beta \in \Zed \rangle$$
given by mapping $f \longmapsto \frac{1}{m} \overline{ev_3}(mf)$. 
\begin{proof}
The relations are given in the comments preceding the Proposition.  Relations (1) and (4) kill
$[t_2^\alpha w_{23}, t_2^\beta w_{23}]$ and $[t_1^\alpha w_{12}, t_1^\beta w_{12}]$ respectively. 
Using Relations (1) and (4) we can simplify relations (2) and (3) to 3-term relations, both expressing
$[t_1^\alpha w_{13}, t_1^\beta w_{13}]$ in the $\Zed[t_1^\pm, t_2^\pm]$-linear span of
$[w_{12}, w_{23}]$. Comparing the two gives the relation
$$\left(t_1^{\alpha-\beta}t_3^{-\beta} - t_1^\alpha t_3^{\alpha-\beta} + (-1)^{n} \left( t_1^\beta t_3^{\beta-\alpha} - t_1^{\beta-\alpha} t_3^{-\alpha} \right)\right) [w_{12},w_{23}] = 0.$$ 

\end{proof}
\end{proposition}

\begin{figure}[H]
{
\psfrag{1}[tl][tl][0.7][0]{$1$}
\psfrag{4}[tl][tl][0.7][0]{$4$}
\psfrag{8}[tl][tl][0.7][0]{$8$}
\psfrag{12}[tl][tl][0.7][0]{$12$}
\psfrag{+}[tl][tl][0.7][0]{$+$}
\psfrag{-}[tl][tl][0.7][0]{$-$}
\psfrag{mn}[tl][tl][0.7][0]{$(-1)^{n-1}$}
\psfrag{mn1}[tl][tl][0.7][0]{$(-1)^{n}$}
\psfrag{t1}[tl][tl][0.7][0]{$t_1$}
\psfrag{t2}[tl][tl][0.7][0]{$t_3$}
\psfrag{t2mk}[tl][tl][0.7][0]{$t^{2-k}$}
\psfrag{t1mk}[tl][tl][0.7][0]{$t^{1-k}$}
\psfrag{tkm1}[tl][tl][0.7][0]{$t^{k-1}$}
\psfrag{tkm2}[tl][tl][0.7][0]{$t^{k-2}$}
\psfrag{BBn}[tl][tl][0.7][0]{$\textcolor{blue}{D^{n-3}}$}
\psfrag{RBn}[tl][tl][0.7][0]{$\textcolor{red}{D^{n-3}}$}
\psfrag{c2i}[tl][tl][0.7][0]{$C_2'[I]$}
\psfrag{c3i}[tl][tl][0.7][0]{$C_3[I]$}
\psfrag{mon1}[tl][tl][0.7][0]{$-t_1^{2-k} t_3^1$}
\psfrag{mon2}[tl][tl][0.7][0]{$(-1)^{n} t_1^{1-k} t_3^{-1}$}
\psfrag{mon3}[tl][tl][0.7][0]{$0$}
\psfrag{mon4}[tl][tl][0.7][0]{$0$}
\psfrag{mon5}[tl][tl][0.7][0]{$+t_1^{-1}t_3^{1-k}$}
\psfrag{mon6}[tl][tl][0.7][0]{$(-1)^{n-1} t_1^1 t_3^{2-k}$}
\psfrag{lab1}[tl][tl][0.7][0]{$lk(\textcolor{red}{t^\alpha Co_2^1}, \textcolor{blue}{t^{\beta-\alpha} Co_1^3})$}
\psfrag{lab2}[tl][tl][0.7][0]{$lk(\textcolor{blue}{t^\alpha Co_2^1}, \textcolor{red}{t^{\beta-\alpha} Co_1^3})$}
\psfrag{lab3}[tl][tl][0.7][0]{$-lk(\textcolor{red}{t^\alpha Co_2^1}, \textcolor{blue}{t^\beta Co_2^3})$}
\psfrag{lab4}[tl][tl][0.7][0]{$-lk(\textcolor{blue}{t^\alpha Co_2^1}, \textcolor{red}{t^\beta Co_2^3})$}
\psfrag{lab5}[tl][tl][0.7][0]{$lk(\textcolor{red}{t^{\alpha-\beta} Co_3^1}, \textcolor{blue}{t^\beta Co_2^3})$}
\psfrag{lab6}[tl][tl][0.7][0]{$lk(\textcolor{blue}{t^{\alpha-\beta} Co_3^1}, \textcolor{red}{t^\beta Co_2^3})$}
$$\includegraphics[width=16cm]{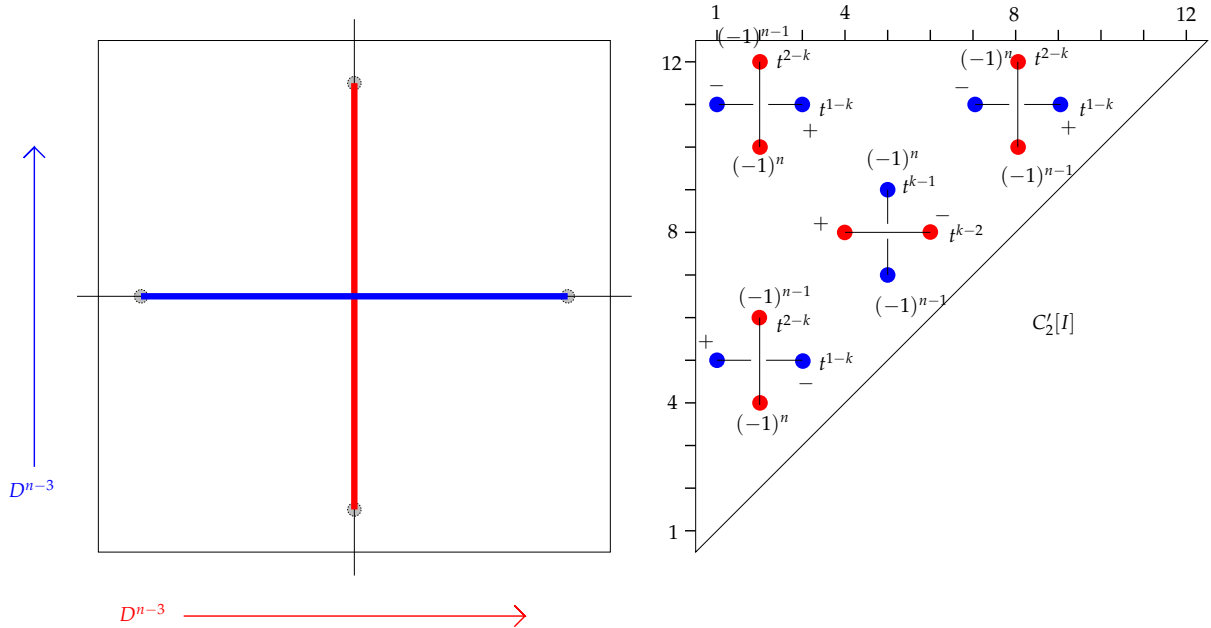}$$
\caption{\label{doublepts1} Scanning through the first $(k-2)$ green dots.}
}
\end{figure}

To compute the $W_3$ invariant, we first consider the homotopy-class of the map on the
$2^{nd}$ stage, $\overline{ev_2}(\delta_k)$. In our family, depicted in Figure \ref{twelvedots}, 
there are $(k-1)$ green dots where double-points are available.  The first $(k-2)$ produce identical
double-point data and they are depicted in Figure \ref{doublepts1}.

In Figure \ref{doublepts1} we have depicted the pre-image of the cohorizontal manifolds for the $2^{nd}$-stage
map, expressed as a map of the form $D^{n-3} \times D^{n-3} \times C_2'[I] \to C_2'[S^1 \times D^{n-1}]$.  While $C_2'[I]$ is technically
a hexagon, the cohorizontal manifold is disjoint from the boundary so for the purpose of exposition we have collapsed
$C_2'[I]$ down to a triangle $\Delta^2$.  The cohorizontal manifolds
are spheres of dimension $n-3$, having a natural surgery `product' decomposition 
$S^{n-3} \equiv D^{n-3} \times \partial I \cup S^{n-4} \times I$.
In the figure, this surgery decomposition is represented by the solid colored arcs on the left side of the figure, depicting
a copy of $D^{n-3}$, together with the pair of similarly-coloured points on the right side of the figure -- for the
$D^{n-3} \times \partial I$ portion.  For the $I \times S^{n-4}$ portion, the spherical boundary is depicted by a pair of large
gray dots on the left-side of the picture, while the black intervals on the right-side of the picture describe an interval $I$.
The $S^{n-4} \times I$ factors occur during the `end homotopy' in the construction of our family (depicted in the top-right
portion of Figure \ref{bbscan}), while the $D^{n-3} \times \partial I$
portion come from the double-points that persist when varying the arcs in the opposite cuff, i.e. depicted in the bottom-right
portion of Figure \ref{bbscan}.

We deduce Figure \ref{doublepts1} from Figure \ref{twelvedots}. The key idea is there are only cohorizontal points (i.e. double points)
for small families where the scanning arc passes through the $k-1$ green dots.  These correspond to the centres of the copies of 
$D^{n-3}$ in the $D^{n-3} \times D^{n-3}$-family corresponding to grabbing the two strands of the embedding in the red or blue 
cylinder respectively, and sweeping them around the barbell and over the respective red or blue cuff. Thus the double-points occur
when the two strands in the coloured cylinders over or undercross the strands running through the bar or the embedded barbell.  The 
numbers decorating features of the lower left (and right) part of Figure \ref{twelvedots} marks the rough parameter-times where 
cohorizontal points occur.  We use this numbering system in Figure \ref{doublepts1}, i.e. coordinate $(2,4)$ has a red dot decorating
it, meaning it records the cohorizontal points from the first $k-2$ green dots, where the strand decorated by $2$ sweeps around
the point on the embedding decorated by $4$. 

\begin{figure}[H]
\begin{multicols}{2}
\begin{itemize}
\item [\textcolor{blue}{(1,5)}] $+t_1 \wedge t_2 \wedge b_1 \wedge \cdots \wedge b_{n-3}$
\item [\textcolor{blue}{(3,5)}] $-t_1 \wedge t_2 \wedge b_1 \wedge \cdots \wedge b_{n-3}$
\item [\textcolor{blue}{(1,11)}] $-t_1 \wedge t_2 \wedge b_1 \wedge \cdots \wedge b_{n-3}$
\item [\textcolor{blue}{(3,11)}] $+t_1 \wedge t_2 \wedge b_1 \wedge \cdots \wedge b_{n-3}$
\item [\textcolor{blue}{(7,11)}] $-t_1 \wedge t_2 \wedge b_1 \wedge \cdots \wedge b_{n-3}$
\item [\textcolor{blue}{(9,11)}] $+t_1 \wedge t_2 \wedge b_1 \wedge \cdots \wedge b_{n-3}$
\item [\textcolor{blue}{(5,7)}] $(-1)^{n-1} t_1 \wedge t_2 \wedge b_1 \wedge \cdots \wedge b_{n-3}$
\item [\textcolor{blue}{(5,9)}] $(-1)^{n} t_1 \wedge t_2 \wedge b_1 \wedge \cdots \wedge b_{n-3}$
\end{itemize}
\columnbreak
\begin{itemize}
\item [\textcolor{red}{(2,4)}] $(-1)^{n} t_1 \wedge t_2 \wedge r_1 \wedge \cdots \wedge r_{n-3}$
\item [\textcolor{red}{(2,6)}] $(-1)^{n-1} t_1 \wedge t_2 \wedge r_1 \wedge \cdots \wedge r_{n-3}$
\item [\textcolor{red}{(2,10)}] $(-1)^{n} t_1 \wedge t_2 \wedge r_1 \wedge \cdots \wedge r_{n-3}$
\item [\textcolor{red}{(2,12)}] $(-1)^{n-1} t_1 \wedge t_2 \wedge r_1 \wedge \cdots \wedge r_{n-3}$
\item [\textcolor{red}{(8,10)}] $(-1)^{n-1} t_1 \wedge t_2 \wedge r_1 \wedge \cdots \wedge r_{n-3}$
\item [\textcolor{red}{(8,12)}] $(-1)^{n} t_1 \wedge t_2 \wedge r_1 \wedge \cdots \wedge r_{n-3}$
\item [\textcolor{red}{(4,8)}] $+ t_1 \wedge t_2 \wedge r_1 \wedge \cdots \wedge r_{n-3}$
\item [\textcolor{red}{(6,8)}] $- t_1 \wedge t_2 \wedge r_1 \wedge \cdots \wedge r_{n-3}$
\end{itemize}
\end{multicols}
\caption{\label{doublepts1s} Normal orientations for the cohorizontal manifolds in the first $(k-2)$ green dots.}
\end{figure}

Notice that the spheres in Figure \ref{doublepts1} are {\bf unlinked} and trivially framed, meaning the $2^{nd}$-stage map is null
homotopic.  Given that the spheres are essentially linearly embedded, you
could think about how the spanning disc for one sphere intersects the others, and some of these intersections are non-trivial.  That said, 
they are avoidable.  Figure \ref{doublepts1} also has some sign information recorded.  This describes the orientations inherited by 
the spheres.  For example, take one of the `red' spheres.  In its surgery decomposition, it consists of the red copy of $D^{n-3} \times \partial I$
together with a copy of $S^{n-4} \times I$.  The sign information in the figure indicates the orientation of the $D^{n-3} \times \partial I$
components.  Since the $D^{n-3} \times \partial I$ components are parallel to the blue $D^{n-3}$ coordinate axis, the sign is a reference
to the orientation of the normal bundle of that portion of the manifold, given in reference to the standard orientation of 
$\Delta^2 \times D^{n-3} \times D^{n-3}$ using the coordinates $(r_1,\cdots, r_{n-3}, b_1, \cdots, b_{n-3}, t_1, t_2)$ in that order, i.e.
and we use the orientation form
$$ dt_1 \wedge dt_2 \wedge dr_1 \wedge \cdots \wedge dr_{n-3} \wedge db_1 \wedge \cdots \wedge db_{n-3}$$

\begin{figure}[H]
{
\psfrag{+}[tl][tl][0.7][0]{$+$}
\psfrag{-}[tl][tl][0.7][0]{$-$}
\psfrag{mn}[tl][tl][0.7][0]{$(-1)^{n-1}$}
\psfrag{mn1}[tl][tl][0.7][0]{$(-1)^{n}$}
\psfrag{t1}[tl][tl][0.7][0]{$t_1$}
\psfrag{t2}[tl][tl][0.7][0]{$t_3$}
\psfrag{t2mk}[tl][tl][0.7][0]{$t^{2-k}$}
\psfrag{t1mk}[tl][tl][0.7][0]{$t^{1-k}$}
\psfrag{tkm1}[tl][tl][0.7][0]{$t^{k-1}$}
\psfrag{tkm2}[tl][tl][0.7][0]{$t^{k-2}$}
\psfrag{BBn}[tl][tl][0.7][0]{$\textcolor{blue}{D^{n-3}}$}
\psfrag{RBn}[tl][tl][0.7][0]{$\textcolor{red}{D^{n-3}}$}
\psfrag{c2i}[tl][tl][0.7][0]{$C_2[I]$}
\psfrag{c3i}[tl][tl][0.7][0]{$C_3[I]$}
\psfrag{mon1}[tl][tl][0.7][0]{$(-1)^{n-1}t_1^{2-k} t_3^{1-k}$}
\psfrag{mon2}[tl][tl][0.7][0]{$-t_1^{1-k} t_3^{2-k}$}
\psfrag{mon3}[tl][tl][0.7][0]{$-t_1^{2-k} t_3$}
\psfrag{mon4}[tl][tl][0.7][0]{$(-1)^{n}\left( t_1^{1-k}t_3^{-1} + t_1^{k-1} t_3 \right)$}
\psfrag{mon5}[tl][tl][0.7][0]{$+t_1 t_3^{k-1} + t_1^{-1} t_3^{1-k} $}
\psfrag{mon6}[tl][tl][0.7][0]{$(-1)^{n-1} t_1 t_3^{2-k} $}
\psfrag{lab1}[tl][tl][0.7][0]{$-lk(\textcolor{red}{t^\alpha Co_2^1}, \textcolor{blue}{t^\beta Co_2^3})$}
\psfrag{lab2}[tl][tl][0.7][0]{$-lk(\textcolor{blue}{t^\alpha Co_2^1}, \textcolor{red}{t^\beta Co_2^3})$}
\psfrag{lab3}[tl][tl][0.7][0]{$lk(\textcolor{red}{t^\alpha Co_2^1}, \textcolor{blue}{t^{\beta-\alpha} Co_1^3})$}
\psfrag{lab4}[tl][tl][0.7][0]{$lk(\textcolor{blue}{t^\alpha Co_2^1}, \textcolor{red}{t^{\beta-\alpha} Co_1^3})$}
\psfrag{lab5}[tl][tl][0.7][0]{$lk(\textcolor{red}{t^{\alpha-\beta} Co_3^1}, \textcolor{blue}{t^\beta Co_2^3})$}
\psfrag{lab6}[tl][tl][0.7][0]{$lk(\textcolor{blue}{t^{\alpha -\beta} Co_3^1}, \textcolor{red}{t^\beta Co_2^3})$}
$$\includegraphics[width=16cm]{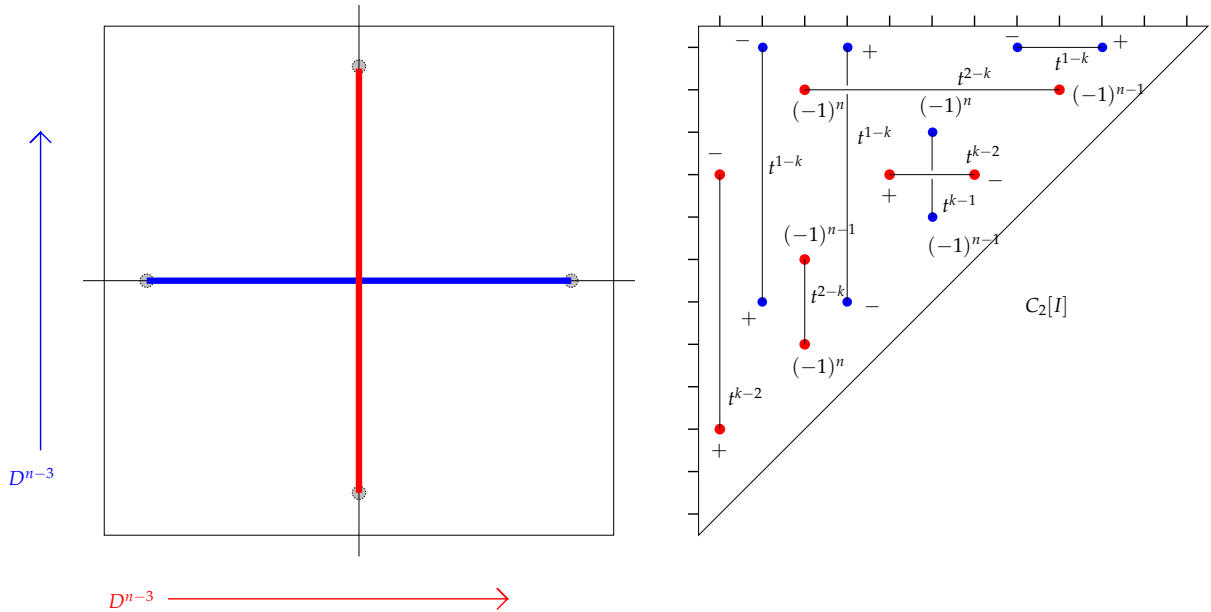}$$
\caption{\label{doublepts2} Scanning through the last green dot.}
}
\end{figure}

With these conventions, the normal orientations for the first $1 \leq l < k-1$ terms are given in Figure \ref{doublepts1s}.  We abbreviate the
$1$-forms $dx$ simply by the symbol $x$. For the last green dot, we have a somewhat different $2^{nd}$-stage diagram given in Figure \ref{doublepts2}.

In Figure \ref{doublepts1s} the point $(1,5)$ and $(3,5)$ are labelled in blue, with 
signs $+$ and $-$ respectively.  In Figure \ref{doublepts1} they are connected by an arc decorated with the monomial 
$t^{1-k}$.  This means that they are part of the preimage of the $t^{1-k} Co_1^2$ manifold.  The plus sign indicates the normal 
orientation of this manifold is $+r_1 \wedge \cdots \wedge r_{n-3}$, i.e. agreeing with the orientation induced by the natural ordering
of the coordinates listed in the order $(t_1,t_2,r_1, \cdots, r_{n-3}, b_1, \cdots, b_{n-3})$. 
Similarly, the disc corresponding to the red dot at $(4,8)$ is labelled in the preimage of $t^{k-2} Co_1^2$ with orientation
described in the Figure \ref{doublepts1s} table.

\begin{figure}[H]
\begin{multicols}{2}
\begin{itemize}
\item [\textcolor{blue}{(2,6)}] $+t_1 \wedge t_2 \wedge b_1 \wedge \cdots \wedge b_{n-3}$
\item [\textcolor{blue}{(2,12)}] $-t_1 \wedge t_2 \wedge b_1 \wedge \cdots \wedge b_{n-3}$
\item [\textcolor{blue}{(4,6)}] $-t_1 \wedge t_2 \wedge b_1 \wedge \cdots \wedge b_{n-3}$
\item [\textcolor{blue}{(4,12)}] $+t_1 \wedge t_2 \wedge b_1 \wedge \cdots \wedge b_{n-3}$
\item [\textcolor{blue}{(8,12)}] $-t_1 \wedge t_2 \wedge b_1 \wedge \cdots \wedge b_{n-3}$
\item [\textcolor{blue}{(10,12)}] $+t_1 \wedge t_2 \wedge b_1 \wedge \cdots \wedge b_{n-3}$
\item [\textcolor{blue}{(6,8)}] $(-1)^{n-1} t_1 \wedge t_2 \wedge b_1 \wedge \cdots \wedge b_{n-3}$
\item [\textcolor{blue}{(6,10)}] $(-1)^{n} t_1 \wedge t_2 \wedge b_1 \wedge \cdots \wedge b_{n-3}$
\end{itemize}
\columnbreak
\begin{itemize}
\item [\textcolor{red}{(3,5)}] $(-1)^{n} t_1 \wedge t_2 \wedge r_1 \wedge \cdots \wedge r_{n-3}$
\item [\textcolor{red}{(3,7)}] $(-1)^{n-1} t_1 \wedge t_2 \wedge r_1 \wedge \cdots \wedge r_{n-3}$
\item [\textcolor{red}{(3,11)}] $(-1)^{n} t_1 \wedge t_2 \wedge r_1 \wedge \cdots \wedge r_{n-3}$
\item [\textcolor{red}{(9,11)}] $(-1)^{n-1} t_1 \wedge t_2 \wedge r_1 \wedge \cdots \wedge r_{n-3}$
\item [\textcolor{red}{(1,3)}] $+t_1 \wedge t_2 \wedge r_1 \wedge \cdots \wedge r_{n-3}$
\item [\textcolor{red}{(1,9)}] $-t_1 \wedge t_2 \wedge r_1 \wedge \cdots \wedge r_{n-3}$
\item [\textcolor{red}{(5,9)}] $+t_1 \wedge t_2 \wedge r_1 \wedge \cdots \wedge r_{n-3}$
\item [\textcolor{red}{(7,9)}] $-t_1 \wedge t_2 \wedge r_1 \wedge \cdots \wedge r_{n-3}$
\end{itemize}
\end{multicols}
\caption{\label{doublepts2s} Normal orientations for the cohorizontal manifolds in the last green dot.}
\end{figure}

\begin{figure}[H]
{
\psfrag{+}[tl][tl][0.7][0]{$+$}
\psfrag{-}[tl][tl][0.7][0]{$-$}
\psfrag{mn}[tl][tl][0.7][0]{$(-1)^{n-1}$}
\psfrag{mn1}[tl][tl][0.7][0]{$(-1)^{n}$}
\psfrag{t1}[tl][tl][0.7][0]{$t_1$}
\psfrag{t2}[tl][tl][0.7][0]{$t_3$}
\psfrag{t2mk}[tl][tl][0.7][0]{$t^{2-k}$}
\psfrag{t1mk}[tl][tl][0.7][0]{$t^{1-k}$}
\psfrag{tkm1}[tl][tl][0.7][0]{$t^{k-1}$}
\psfrag{tkm2}[tl][tl][0.7][0]{$t^{k-2}$}
\psfrag{BBn}[tl][tl][0.7][0]{$\textcolor{blue}{D^{n-3}}$}
\psfrag{RBn}[tl][tl][0.7][0]{$\textcolor{red}{D^{n-3}}$}
\psfrag{c2i}[tl][tl][0.7][0]{$C_2[I]$}
\psfrag{c3i}[tl][tl][0.7][0]{$C_3[I]$}
\psfrag{mon1}[tl][tl][0.7][0]{$-t_1^{2-k} t_3^1$}
\psfrag{mon2}[tl][tl][0.7][0]{$(-1)^{n} t_1^{1-k} t_3^{-1}$}
\psfrag{mon3}[tl][tl][0.7][0]{$0$}
\psfrag{mon4}[tl][tl][0.7][0]{$0$}
\psfrag{mon5}[tl][tl][0.7][0]{$+t_1^{-1}t_3^{1-k}$}
\psfrag{mon6}[tl][tl][0.7][0]{$(-1)^{n-1} t_1^1 t_3^{2-k}$}
\psfrag{lab1}[tl][tl][0.7][0]{$lk(\textcolor{red}{t^\alpha Co_2^1}, \textcolor{blue}{t^{\beta-\alpha} Co_1^3})$}
\psfrag{lab2}[tl][tl][0.7][0]{$lk(\textcolor{blue}{t^\alpha Co_2^1}, \textcolor{red}{t^{\beta-\alpha} Co_1^3})$}
\psfrag{lab3}[tl][tl][0.7][0]{$-lk(\textcolor{red}{t^\alpha Co_2^1}, \textcolor{blue}{t^\beta Co_2^3})$}
\psfrag{lab4}[tl][tl][0.7][0]{$-lk(\textcolor{blue}{t^\alpha Co_2^1}, \textcolor{red}{t^\beta Co_2^3})$}
\psfrag{lab5}[tl][tl][0.7][0]{$lk(\textcolor{red}{t^{\alpha-\beta} Co_3^1}, \textcolor{blue}{t^\beta Co_2^3})$}
\psfrag{lab6}[tl][tl][0.7][0]{$lk(\textcolor{blue}{t^{\alpha-\beta} Co_3^1}, \textcolor{red}{t^\beta Co_2^3})$}
$$\includegraphics[width=14cm]{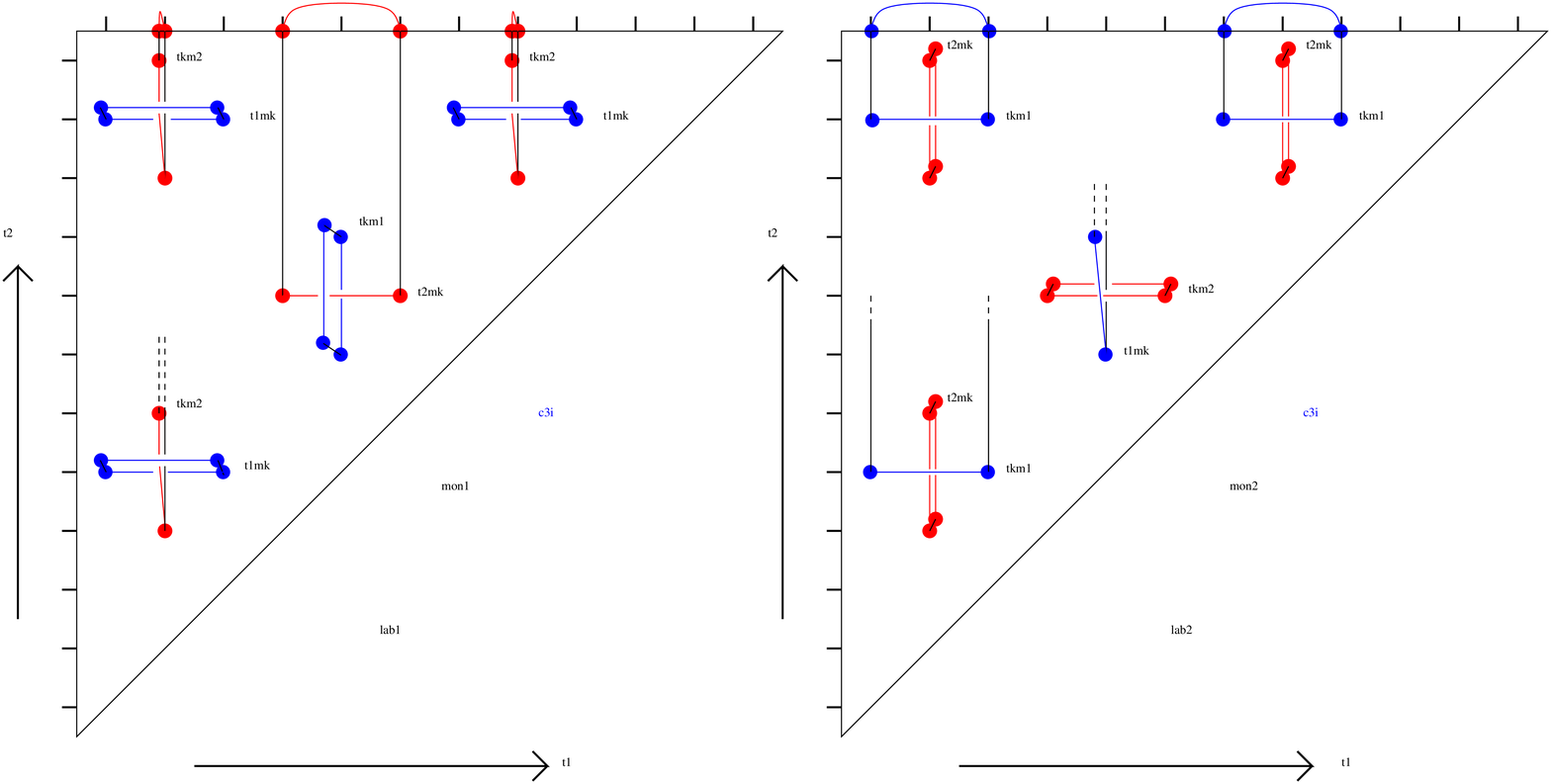}$$
$$\includegraphics[width=14cm]{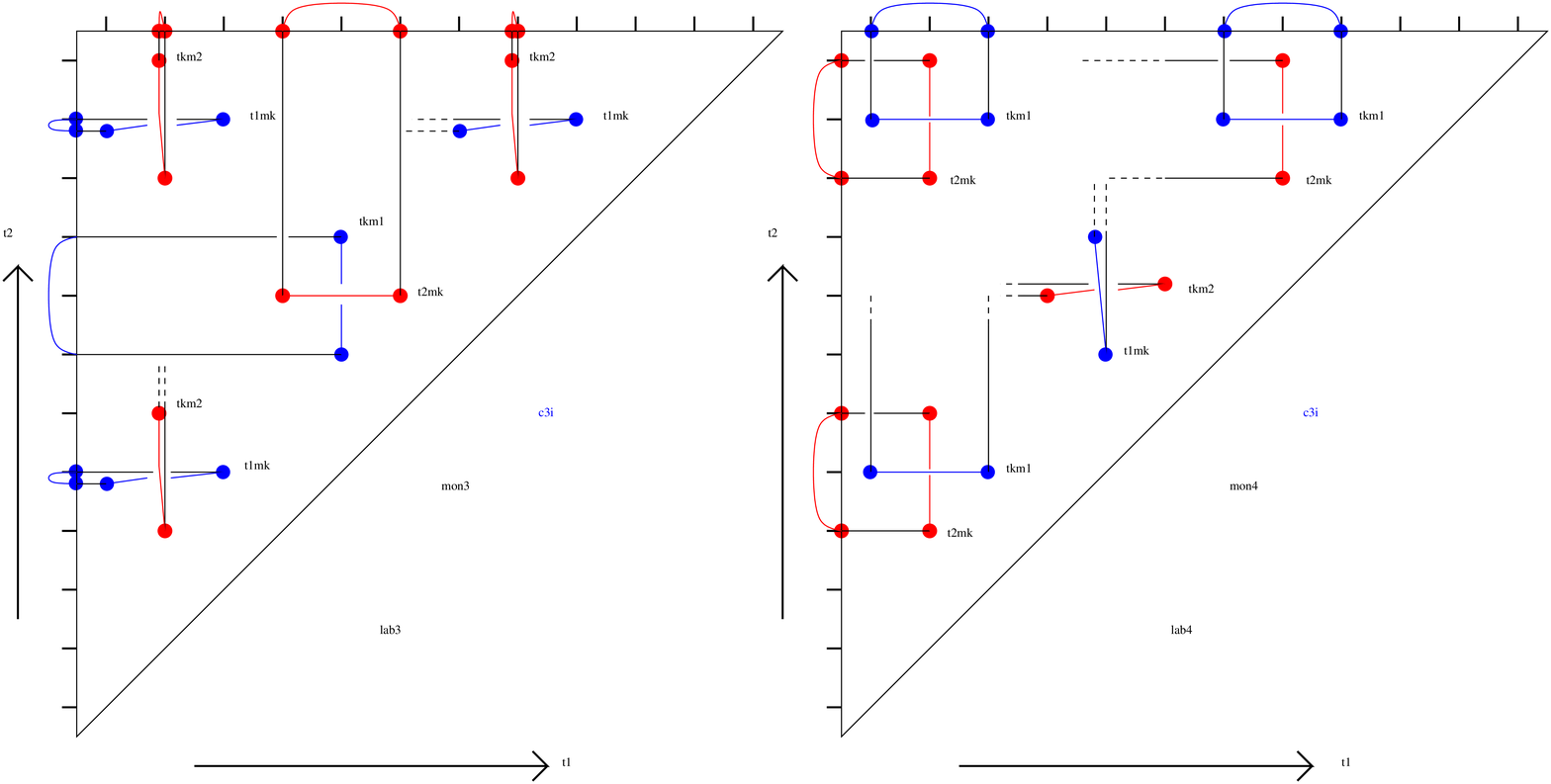}$$
$$\includegraphics[width=14cm]{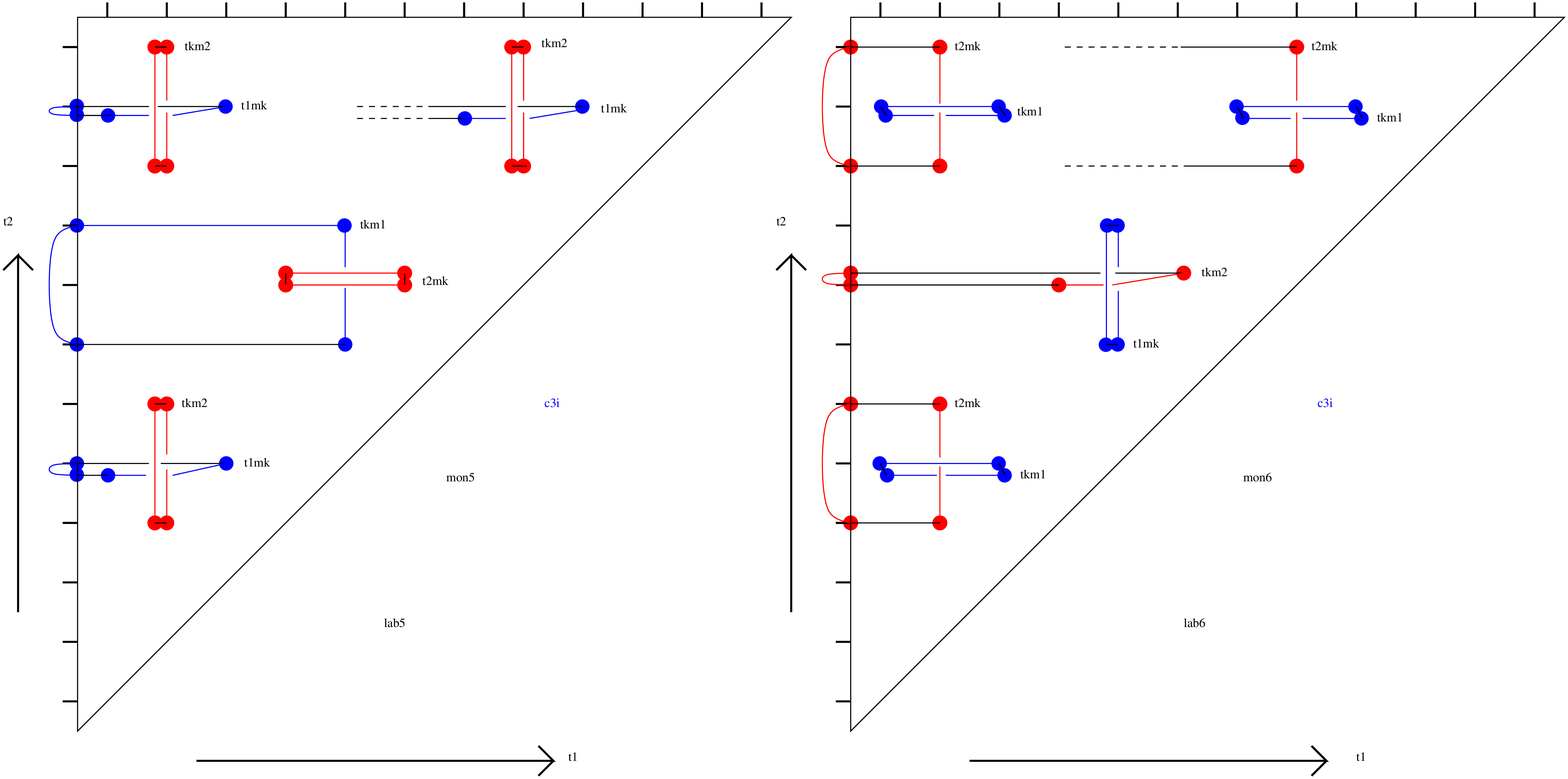}$$
}
\caption{\label{cohfirst}Cohorizontal manifolds for first $k-2$ green dots, in $D^{n-3} \times D^{n-3} \times C_3'[I]$}
\end{figure}

\begin{figure}[H]
{
\psfrag{+}[tl][tl][0.7][0]{$+$}
\psfrag{-}[tl][tl][0.7][0]{$-$}
\psfrag{mn}[tl][tl][0.7][0]{$(-1)^{n-1}$}
\psfrag{mn1}[tl][tl][0.7][0]{$(-1)^{n}$}
\psfrag{t1}[tl][tl][0.7][0]{$t_1$}
\psfrag{t2}[tl][tl][0.7][0]{$t_3$}
\psfrag{t2mk}[tl][tl][0.7][0]{$t^{2-k}$}
\psfrag{t1mk}[tl][tl][0.7][0]{$t^{1-k}$}
\psfrag{tkm1}[tl][tl][0.7][0]{$t^{k-1}$}
\psfrag{tkm2}[tl][tl][0.7][0]{$t^{k-2}$}
\psfrag{BBn}[tl][tl][0.7][0]{$\textcolor{blue}{D^{n-3}}$}
\psfrag{RBn}[tl][tl][0.7][0]{$\textcolor{red}{D^{n-3}}$}
\psfrag{c2i}[tl][tl][0.7][0]{$C_2[I]$}
\psfrag{c3i}[tl][tl][0.7][0]{$C_3[I]$}
\psfrag{mon1}[tl][tl][0.7][0]{$(-1)^{n-1}t_1^{2-k} t_3^{1-k}$}
\psfrag{mon2}[tl][tl][0.7][0]{$-t_1^{1-k} t_3^{2-k}$}
\psfrag{mon3}[tl][tl][0.7][0]{$-t_1^{2-k} t_3$}
\psfrag{mon4}[tl][tl][0.7][0]{$(-1)^{n}\left( t_1^{1-k}t_3^{-1} + t_1^{k-1} t_3 \right)$}
\psfrag{mon5}[tl][tl][0.7][0]{$+t_1 t_3^{k-1} + t_1^{-1} t_3^{1-k} $}
\psfrag{mon6}[tl][tl][0.7][0]{$(-1)^{n-1} t_1 t_3^{2-k} $}
\psfrag{lab1}[tl][tl][0.7][0]{$-lk(\textcolor{red}{t^\alpha Co_2^1}, \textcolor{blue}{t^\beta Co_2^3})$}
\psfrag{lab2}[tl][tl][0.7][0]{$-lk(\textcolor{blue}{t^\alpha Co_2^1}, \textcolor{red}{t^\beta Co_2^3})$}
\psfrag{lab3}[tl][tl][0.7][0]{$lk(\textcolor{red}{t^\alpha Co_2^1}, \textcolor{blue}{t^{\beta-\alpha} Co_1^3})$}
\psfrag{lab4}[tl][tl][0.7][0]{$lk(\textcolor{blue}{t^\alpha Co_2^1}, \textcolor{red}{t^{\beta-\alpha} Co_1^3})$}
\psfrag{lab5}[tl][tl][0.7][0]{$lk(\textcolor{red}{t^{\alpha-\beta} Co_3^1}, \textcolor{blue}{t^\beta Co_2^3})$}
\psfrag{lab6}[tl][tl][0.7][0]{$lk(\textcolor{blue}{t^{\alpha -\beta} Co_3^1}, \textcolor{red}{t^\beta Co_2^3})$}
$$\includegraphics[width=14cm]{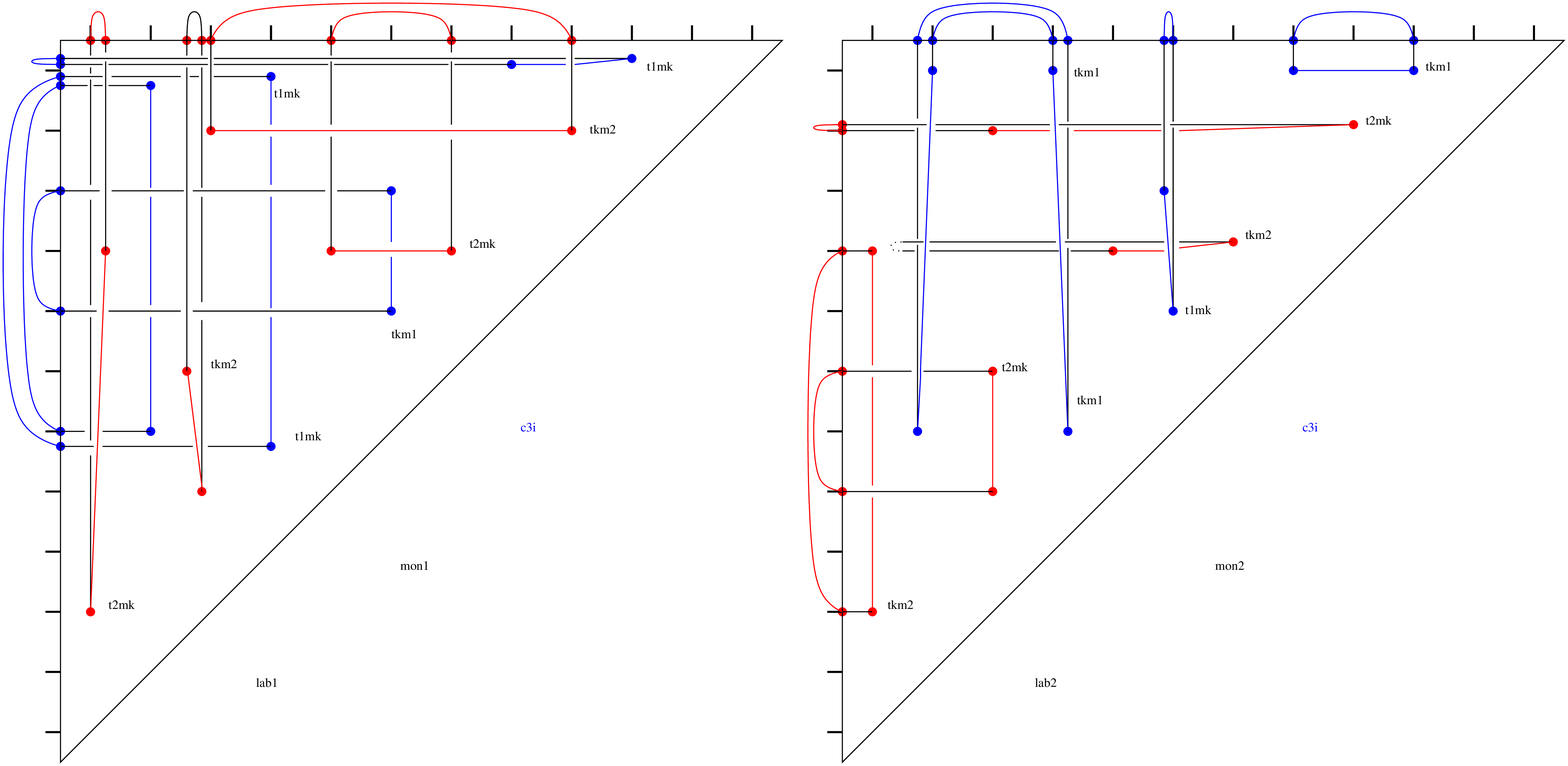}$$
$$\includegraphics[width=14cm]{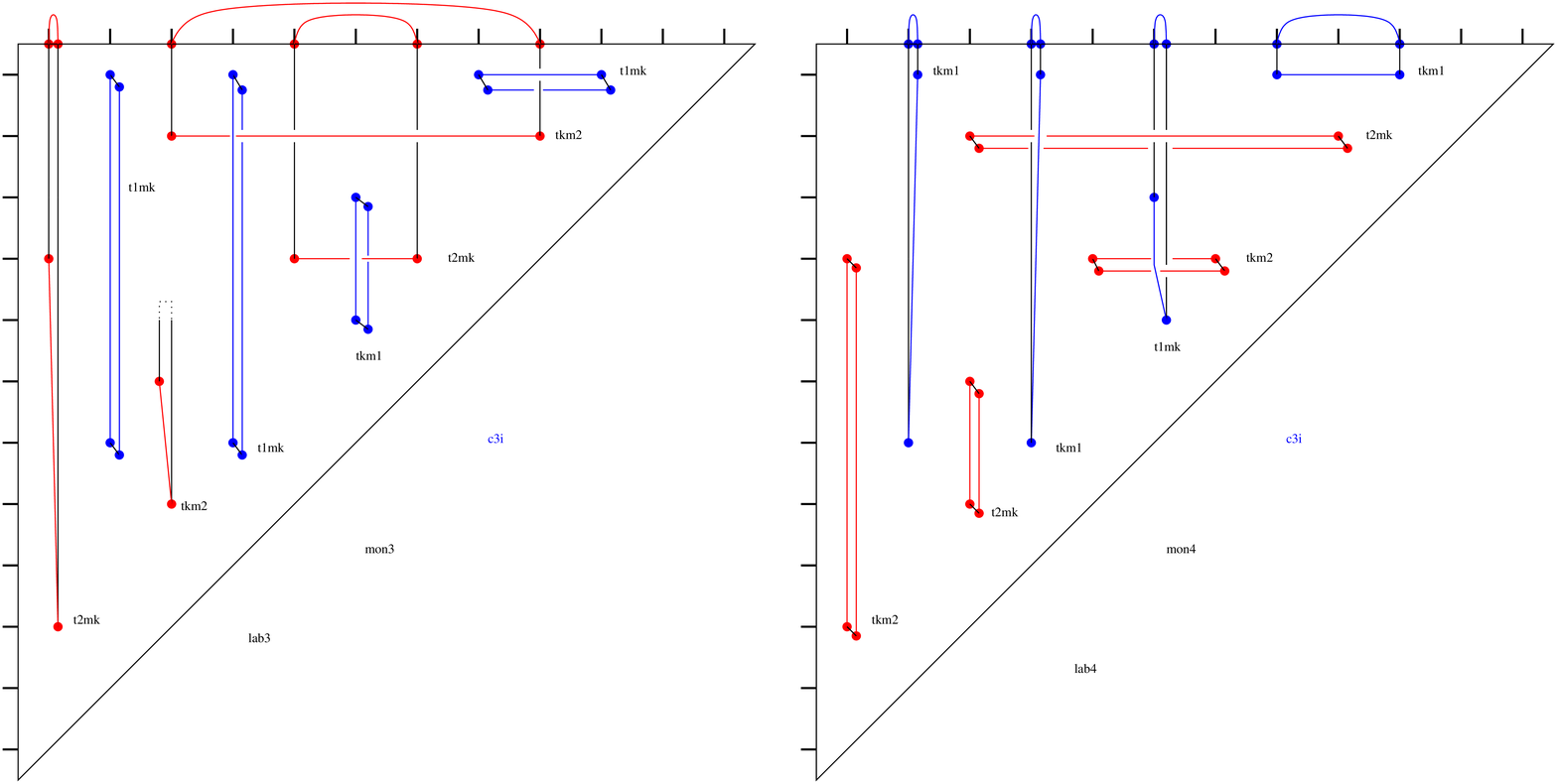}$$
$$\includegraphics[width=14cm]{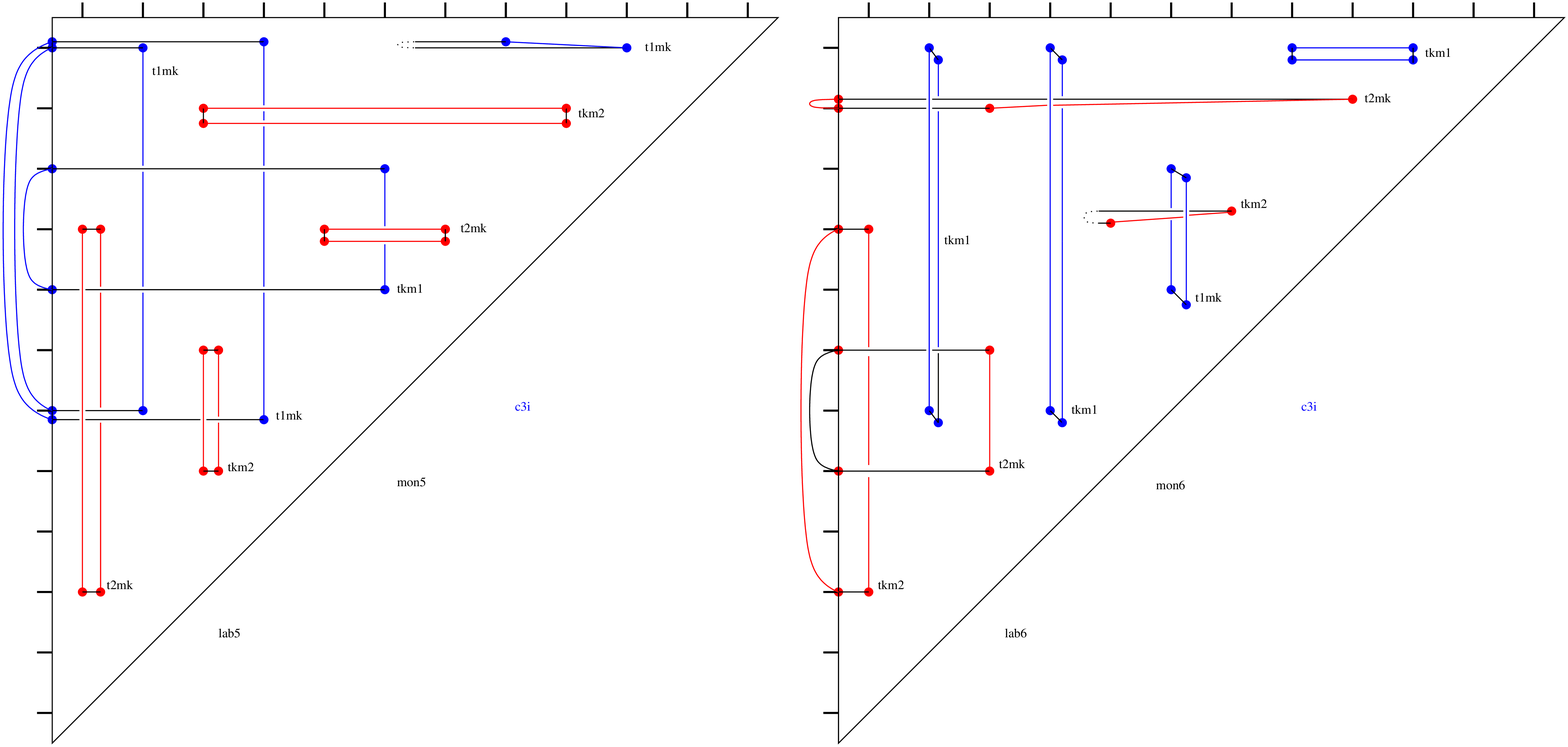}$$
}
\caption{\label{cohsecond}Cohorizontal manifolds for the last green dot, in $D^{n-3} \times D^{n-3} \times C_3'[I]$}
\end{figure}

We construct the $3^{rd}$-stage map cohorizontal manifolds from the $2^{nd}$-stage map, the idea being that
whichever cohorizontal manifold one is considering, it will be constant in one of the three parameters of
$C_3'[I]$, given that the cohorizontal condition is a constraint on only two of the three coordinates of
$C_3'[I]$.  Attaching the null-homotopies described for the $2^{nd}$-stage map allows us to close-off the
cohorizontal manifolds, getting a collection of disjoint spheres in 
a neighbourhood of $D^{n-3} \times D^{n-3} \times C_3'[I]$. We emphasize `neighbourhood' since the null-homotopy
attachments are external to $D^{n-3} \times D^{n-3} \times C_3'[I]$.  In Figures \ref{cohfirst} and \ref{cohsecond} we suppress the
$D^{n-3} \times D^{n-3}$ factors, since it would only be a repetition of Figure \ref{doublepts1}. In our diagrams
we display the $t_1$ and $t_3$ coordinates, with $t_2$ being out of the page.  Thus our diagram depicts collections
of spheres diffeomorphic to $S^{n-2}$.  We simplify our sketch of $C_3'[I]$ to be a simple tetrahedron $\Delta^3$.  
If one looks at the diagram, one sees a collection of disjoint circles in a neighbourhood of $\Delta^3$, some linking,
and others not.  As with Figures \ref{doublepts1} and \ref{doublepts2}, these diagrams are in `product form', i.e. 
these spheres have the form $D^{n-3} \times S^1 \cup S^{n-4} \times D^2$, where the circle and $D^2$ factors live
in the neighbourhood of $\Delta^3$, and the $D^{n-3}$ factor comes from one of the factors in the parametrization of our family.  
Thus our figures only depict the $\{0\} \times S^1$ portions of our spheres, but fortunately for us, this is precisely
where our double-points occur.  

The above kind of geometry occurs in the study of standard linking pairs. While this is an elementary geometric observation, 
the first author learned about this phenomenon from Haefliger \cite{Ha}.  For example, if we take a standard
linking pair in $S^{n-1}$ with $i+j=n$, i.e. 
$$S^{n-1} \equiv \partial D^n \equiv \partial (D^i \times D^j) = S^{i-1} \times D^j \cup D^i \times S^{j-1}$$
we can go one step further and think of $D^n \times \{0\}$ as the equator in $D^{n+1}$, giving
$$S^n \equiv \partial D^{n+1} \equiv S^i \times D^j \cup D^{i+1} \times S^{j-1} \text{ or } S^{i-1} \times D^{j+1} \cup D^i \times S^j.$$
i.e. one can think of a $(S^{i-1},S^{j-1})$ standard linking pair in $S^{n-1}$ as equatorial in  
a $(S^i, S^{j-1})$ standard linking pair in $S^n$, or the reverse, in a $(S^{i-1}, S^j)$ standard linking pair in $S^n$.  This 
is a single step in an inductive suspension process that can generate all standard linking pairs of spheres, from linking pairs
of the form $(S^0,S^0)$ in $S^1$. 

Figures \ref{cohfirst} and \ref{cohsecond} are depicting standard linking pairs of the form $(S^{n-2}, S^{n-2})$ where we 
are inducting up from the base
case of disjoint linking circles in $\Delta^3$, where for one circle we suspend up using the $D^{n-3} \times \{0\}$ factor, and for the other
circle we suspend up using the $\{0\} \times D^{n-3}$ factor. 

In Figure \ref{cohfirst} we break up the cohorizontal manifolds into their constituent parts, depending on which cuff the points are being
mapped to (via colour) and the translate of the relevant cohorizontal manifold.  In the top-left part of Figure \ref{cohfirst} 
($lk(\textcolor{red}{t^\alpha Co_2^1},\textcolor{blue}{t^{\beta-\alpha} Co_1^3})$) there is only the one linking pair, 
with monomial
$t_1^{2-k} t_3^1$.  To compute the signs, we count the signed overcrossings
of $t^\alpha Co_2^1 - t^{\alpha-\beta} Co_3^1$ over $t^{\beta-\alpha} Co_1^3 - t^\beta Co_2^3$. 
The normal orientation of the arc parallel to the $t_3$-axis through $(6,8,8)$ is
$-t_1 \wedge t_2 \wedge r_1 \wedge \cdots \wedge r_{n-3}$.  The normal orientation to
the arc parallel to the $t_2$-axis through $(5,5,7)$ is $(-1)^{n-1} t_1 \wedge t_3 \wedge b_1 \wedge \cdots \wedge b_{n-3}$. 
Repeating for all six sub-diagrams of Figure \ref{cohfirst}, we get the sum of all these terms for the first $k-2$ green dots as
$$(k-2)\left(t_1^{-1} t_3^{1-k} + (-1)^{n} t_1^{1-k} t_3^{-1} - t_1^{2-k}t_3^1 + (-1)^{n-1} t_1^1 t_3^{2-k}\right).$$

If we repeat for Figure \ref{cohsecond}, the sum of the terms for the last green dot gives
$$ (-1)^{n-1} t_1^{2-k} t_3^{1-k} - 
   t_1^{1-k} t_3^{2-k} - 
   t_1^{2-k} t_3 + 
   (-1)^{n}t_1^{1-k}t_3^{-1} + 
   (-1)^{n} t_1^{k-1} t_3^1 +
   t_1 t_3^{k-1} + 
   t_1^{-1} t_3^{1-k} + 
   (-1)^{n-1} t_1 t_3^{2-k}.$$

Putting this together with the first $k-2$ green dots, we have
$$W_3(\delta_k) = (k-1)\left(t_1^{-1} t_3^{1-k} + (-1)^{n} t_1^{1-k} t_3^{-1} - t_1^{2-k}t_3^1 + (-1)^{n-1} t_1 t_3^{2-k}\right) +$$
 $$ t_1t_3^{k-1} + (-1)^{n} t_1^{k-1} t_3 - t_1^{1-k} t_3^{2-k} + (-1)^{n-1} t_1^{2-k} t_3^{1-k}$$
which completes the proof of Theorem \ref{mainsmooththm}. 

Recall the hexagon relation.

\begin{align}
t_1^{\alpha}t_3^{\beta} + (-1)^{n} t_1^{-\beta} t_3^{-\alpha} & = t_1^{\alpha-\beta} t_3^{\alpha} + (-1)^{n} t_1^{-\alpha} t_3^{\beta-\alpha} \\
  & = t_1^{-\beta} t_3^{\alpha-\beta} + (-1)^{n} t_1^{\beta-\alpha} t_3^{\beta} \\
  & = t_1^{-\alpha} t_3^{-\beta} + (-1)^{n} t_1^\beta t_3^\alpha \\
  & = t_1^{\beta-\alpha} t_3^{-\alpha} + (-1)^{n} t_1^\alpha t_3^{\alpha-\beta} \\
  & = t_1^\beta t_3^{\beta-\alpha} + (-1)^{n} t_1^{\alpha-\beta} t_3^{-\beta}
\end{align}

Plugging in $k=3$ in the above formula for $W_3(\delta_3)$ gives
$$W_3(\delta_3) = \begin{cases} 0 & \text{ if } n \text{ even} \\
  4t_1^{-2}t_3^{-1} + 4t_1^{-1}t_3^{-2}-2t_1^{-1}t_3-2t_1t_3^{-1} & \text{ if } n \text{ odd.} \end{cases}$$

\section{Homeomorphisms of $S^1 \times D^{n-1}$}\label{homeosec}


In this section we give a proof of Theorem \ref{main}.

\noindent{\bf Proof of Theorem \ref{main}:} We supply an argument that a quotient of the $W_3$ invariant is definable out of 
$\pi_{n-4} \Homeo(S^1 \times D^{n-1})$.  This will suffice to show $\pi_{n-4} \Homeo(S^1 \times D^{n-1})$
is not finitely-generated for all $n \geq 4$.  

Let $\Emb^\tau(I, S^1 \times D^{n-1})$ denote the space of topological embeddings of $I = [0,1]$ in $S^1 \times D^{n-1}$. 
We require these embeddings send $0$ to $(1,-*)$ and $1$ to $(1,*)$ where $* \in \partial D^{n-1}$ is a choice of basepoint.
We similarly require that the embedding send the interior of $I$ to the interior of $S^1 \times D^{n-1}$. This last 
condition does not affect the homotopy-type of the space it does make some technical arguments easier to read. We give this
embedding space the compact-open topology.

We let $\Delta^k$ denote the standard simplex, 
$$\Delta^k = \{ (t_1,t_2, \cdots,t_k) \in \Real^k : 0 \leq t_1 \leq t_2 \leq \cdots \leq t_k \leq 1\}.$$

A topological embedding $f : I \to X$ induces a map
$$f_* : \Delta^k \to (S^1 \times D^{n-1})^k$$
defined by $f_*(t_1, \cdots, t_k) = (f(t_1), \cdots, f(t_k))$.  We list its properties.  Given a set
$A = \{i_1,\cdots,i_j\} \subset \{1,2,\cdots,k\}$, the $A$-diagonal of $X^k$ denotes the subspace of 
$X^k$ where $x_{i_1} = x_{i_2} =\cdots = x_{i_j}$.   Call the subspace of $\Delta^k$ satisfying
$t_1=0$ the {\bf initial facet} of $\Delta^k$, and the subspace satisfying $t_k=1$ the {\bf terminal facet}
of $\Delta^k$. The subset of $(S^1 \times D^{n-1})^k$ satisfying $p_1=(1,-*)$ we call the initial 
facet, and $p_k=(1,*)$ the terminal facet of $(S^1 \times D^{n-1})^k$. 

\begin{itemize}
\item[(a)] The induced map $f_*$ sends $A$-diagonals to $A$-diagonals, 
moreover the pre-images of $A$-diagonals are $A$-diagonals. 
\item[(b)] $f_*$ sends the initial facet of $\Delta^k$ to the initial facet of $(S^1 \times D^{n-1})^k$, similarly the
terminal facets. 
\item[(c)] If we lift $f_*$ to a map of universal covers $\Delta^k \to (\Real^1 \times D^{n-1})^k$ then two covering translates
of points of the image agree if and only if the covering translates are identical, i.e. $p_i = t.p_j$ is possible if and only
if $t=0 \in \pi_1 (S^1 \times D^{n-1})^k$ and $p_i=p_j$.
\end{itemize}

Given a manifold $M$ define
$$ C^\tau_k(M) = \{ (p_1,\cdots,p_k) \in \tilde M^k : p_i \notin (\pi_1 M \setminus \{0\}).p_j \ \forall i,j \}$$
where $\tilde M$ is the universal cover of $M$. This definition can be interpreted as saying that
any two listed points $p_i, p_j \in \tilde M$ either have disjoint $\pi_1 M$-orbits, or when the orbits
intersect we have $p_i = p_j$. We call $C^\tau_k(M)$ the {\bf principal} configuration space of $k$ points
in the universal cover of $M$. Thus item (c) above states the lift of $f_* : \Delta^k \to (S^1 \times D^{n-1})^k$ to
the universal cover is a map of the form $\Delta^k \to C_k^\tau(S^1 \times D^{n-1})$. Items (a) and (b) should be thought
of as a relative mapping space condition. 

As a space, $C^\tau_k(M)$ is the orbit configuration space of
the universal cover of $M$ union the $\{i,j\}$-diagonals, i.e. thinking of $C^\tau_k(M)$ as a subspace of
$(\tilde M)^k$ it is the union of $C_k(\tilde M)$ with the $\{i,j\}$-diagonals for all $i \neq j$.   
See for example Fred Cohen's work on orbit configuration spaces \cite{Fred}. For our purposes we need to
know $\pi_{2n-1} \tilde C_2^\tau(S^1 \times D^{n-1})$ and enough of 
$\pi_{2n-3} \tilde C_3(S^1 \times D^{n-1}) \otimes \Rat$, which is the content of the next proposition. 

\begin{prop}\label{forgetker}
The homotopy group $\pi_{n-1} \tilde C_2(S^1 \times D^{n-1})$ is freely generated by the elements
$t_2^a w_{12}$ with $a \in \Zed \setminus \{0\}$. 

The kernel of the map 
$$\pi_{2n-3} \tilde C_3(S^1 \times D^{n-1}) \otimes \Rat \to \pi_{2n-3} C^\tau_3(S^1 \times D^{n-1}) \otimes \Rat$$
trivially intersects the $\Rat$-span of elements of the form $[t_j^a w_{ij}, t_j^b w_{jk}]$ for $a,b \in \Zed \setminus \{0\}$, 
$a \neq b$ and $\{i,j,k\} = \{1,2,3\}$.
\begin{proof}
Concerning $\pi_{n-1} \tilde C_2(S^1 \times D^{n-1})$, the forgetful map $\tilde C_2(S^1 \times D^{n-1}) \to \Real \times D^{n-1}$
is a locally-trivial fiber bundle with fiber diffeomorphic to the complement of the non-trivial covering translates of a point, 
which has the homotopy-type of a wedge of spheres, the generators in dimension $(n-1)$ being our $t_2^a w_{12}$ classes with $a \neq 0$.

Concerning $\pi_{2n-1} \tilde C_2(S^1 \times D^{n-1})$ we consider the collinear manifolds $Col^1_{\alpha, \beta}$ and 
$Col^3_{\alpha, \beta}$ as being subsets of
$C_3^\tau(S^1 \times D^{n-1})$ via the inclusion $C_3^\tau(S^1 \times D^{n-1}) \to (\Real \times D^{n-1})^3$.
The manifolds $Col^1_{\alpha,\beta}$ and $Col^3_{\alpha,\beta}$ are disjoint and closed in $C_3^\tau(S^1 \times D^{n-1})$ 
provided $\alpha \neq 0 \neq \beta$ and $\alpha \neq \beta$. 

If we pull-back the pair $Col^1_{\alpha,\beta}$ and $Col^3_{\alpha,\beta}$ via the map 
$[t_2^a w_{12}, t_2^b w_{23}]$ we get a disjoint oriented manifold pair with linking number 
$\pm 1$ provided $\alpha=-a$ and $\beta=-b$, otherwise we get zero.  Thus the
set of brackets of the form $[t_2^a w_{12}, t_2^b w_{23}]$ are linearly independent in 
$\pi_{2n-3} C^\tau_3(S^1 \times D^{n-1}) \otimes \Rat$ provided $a \neq b$
and $a \neq 0 \neq b$. 
\end{proof}
\end{prop}

In the smooth category, given a diffeomorphism of $S^1 \times D^{n-1}$ we consider the induced scanning map of the disc
$\{1\} \times D^{n-1}$, this gave us an element of $\Omega^{n-2} \Emb(I, S^1 \times D^{n-1})$.  We follow that same outline
for homeomorphisms, using the map $\Homeo(S^1 \times D^{n-1}) \to \Omega^{n-2} \Emb^\tau(I, S^1 \times D^{n-1})$.

In the smooth case, the induced map on the second stage was torsion.  In the topological case, the `second stage' we take
as the map
$$\Emb^\tau(I, S^1 \times D^{n-1}) \to \Map(\Delta^2, C_2^\tau(S^1 \times D^{n-1}))$$
with the associated boundary conditions, i.e. this is a stratum-preserving mapping space, as described in conditions (a) and 
(b). Notice that the forgetful map $C_k^\tau(S^1 \times D^{n-1}) \to C_{k-1}^\tau(S^1 \times D^{n-1})$
is in general not a fibration, but in the case $k=2$ it is, with the fiber having the homotopy-type of 
$\Real \times D^{n-1} \setminus (\Zed \setminus \{0\}).p$ where $p \in int(\Real \times D^{n-1})$, i.e. this has the homotopy-type
of a wedge of spheres.  Thus given an element of $\pi_{n-4} \Homeo(S^1 \times D^{n-1})$ the induced element of the second stage, 
is a stratum-preserving map of the form $$S^{n-4} \times D^{n-2} \times \Delta^2 \to C_2^\tau(S^1 \times D^{n-1}).$$
The restriction of this map to the boundary facets of $\Delta^2$ give canonically null-homotopic maps, thus we can cap-off
the above map to construct a map $S^{2n-4} \to C_2^\tau(S^1 \times D^{n-1})$. Given that $\pi_{2n-4} S^{n-1}$ is torsion when
$n \geq 4$, our map is torsion.  Like in the smooth case, some multiple of the $3^{rd}$-stage map
$$ S^{n-4} \times D^{n-2} \times \Delta^3 \to C^\tau_3(S^1 \times D^{n-1})$$
is null on the boundary.  Will attach a choice of null-homotopy, and as in the smooth case the induced element of
$$\pi_{2n-3} C^\tau_3(S^1 \times D^{n-1})$$
is well-defined up to an error terms coming from a subgroup $R'$. The subgroup $R'$ is the image of $R$ under the 
induced map from the forgetful map
$$\pi_{2n-3} C_3(S^1 \times D^{n-1}) \to \pi_{2n-3} C^\tau(S^1 \times D^{n-1}).$$

So we have a commutative diagram

$$\xymatrix{\pi_{n-4} \Diff(S^1 \times D^{n-1}) \otimes \Rat \ar[d] \ar[r]^-{W_3} & \pi_{2n-3} C_3(S^1 \times D^{n-1}) \otimes \Rat / R \ar[d] \\
 \pi_{n-4} \Homeo(S^1 \times D^{n-1}) \otimes \Rat \ar[r]^-{W'_3} & \pi_{2n-3} C^\tau_3(S^1 \times D^{n-1}) \otimes \Rat / R' }.$$

Due to Proposition \ref{forgetker}, our elements $\delta_k$ satisfy

$$W'_3(\delta_k) = (k-1)\left(t_1^{-1} t_3^{1-k} + (-1)^{n} t_1^{1-k} t_3^{-1} - t_1^{2-k}t_3^1 + (-1)^{n-1} t_1 t_3^{2-k}\right) +$$
 $$ t_1t_3^{k-1} + (-1)^{n} t_1^{k-1} t_3 - t_1^{1-k} t_3^{2-k} + (-1)^{n-1} t_1^{2-k} t_3^{1-k}$$

which are non-trivial and linearly independent for $k \geq 4$.  The key observation is that these elements lie in the $12$-element orbits
of the dihedral group of the hexagon, and these orbits do not belong
to the kernel of the map $\pi_{2n-3} \tilde C_3(S^1 \times D^{n-1}) \otimes \Rat \to \pi_{2n-3} C^\tau_3(S^1 \times D^{n-1}) \otimes \Rat / R'$
by Proposition \ref{forgetker}, completing the proof of Theorem \ref{main}.  \qed

Sander Kupers has informed us one can further prove that the kernel of the map 
$\pi_{2n-6} \Emb(I, S^1 \times D^{n-1}) \to \pi_{2n-6} \Emb^\tau(I, S^1 \times D^{n-1})$ is the image of the
inclusion $\pi_{2n-6} \Emb(I, D^n) \to \pi_{2n-6} \Emb(I, S^1 \times D^{n-1})$. This is the subgroup given by 
embeddings disjoint from $\{-1\} \times D^{n-1}$. This result (unpublished) would give an alternative proof of Theorem \ref{main}.
One can also obtain Theorem \ref{main} in dimensions different from $n=4,5$ and $7$ using smoothing theory \cite{Kupers}. 
More generally, its known that the homotopy groups of $B(\Homeo(M))$ and $B(\Diff(M))$ are finitely generated whenever
$\pi_1 M$ is finite and $M$ has even dimension, different from $4$ \cite{BKK}. 



\section{Automorphisms of complete finite volume hyperbolic manifolds}\label{hypsec}

Farrell and Jones \cite{FaJo} proved that provided $N$ is a compact hyperbolic manifold of dimension greater than or equal to $11$, 
then $\pi_0 \Diff(N)$ and $\pi_0 \Homeo(N)$ are not finitely generated.  In particular, $\Diff(N)$ and $\Homeo(N)$ do not
have the homotopy-type of compact manifolds.  The purpose of this section is to reduce $11$ to $4$.  By  \cite{Ga2}, it 
is known that the diffeomorphism group of a complete finite volume hyperbolic $3$-manifold has the homotopy-type of 
its isometry group, i.e. it has the homotopy-type of a discrete, finite set.  Thus the results of this section are optimal. 
While the Smale Conjecture for hyperbolic 3-manifolds is stated for closed manifolds in the introduction to \cite{Ga2}, it
follows for complete manifolds by Lemma 7.2 and Theorem 7.3 of \cite{Ga2}.

The proof of Farrell and Jones is constrained by two important dimension restrictions: 1) we do not yet know the optimal range for pseudo-isotopy
stability.  Indeed, we still depend on the initial result of Igusa \cite{Igusa}. 2) Farrell and Jones also depend on the work
of Hatcher and Wagoner \cite{HW} which begins in dimension $6$. 

While Farrell and Jones compute the mapping class group of $N$ in the smooth and topological cases, we restrict to
$\pi_{n-4} \Diff(N)$ and $\pi_{n-4} \Homeo(N)$. In a future paper \cite{BG2} we anticipate extending these arguments 
to the level of mapping class groups.

\begin{thm}\label{mainthm}
If $N$ is a complete finite volume  hyperbolic manifold of dimension $n \geq 4$ then both
$$\pi_{n-4} \Diff(N) \text{ and } \pi_{n-4} \Homeo(N)$$ 
are not finitely generated.

\begin{proof}  We give the proof for $n=4$ and $N$ orientable.  
In a neighbourhood $N(\gamma)$ of an embedded  closed  geodesic $\gamma$, implant the  barbell $\delta_k$ to obtain the diffeomorphism $f_k\in \Diff_0(N)$.
  Let $\hat f_k$ be the lift of $f_k$ to the covering space $\hat N_\gamma$ of $N$ corresponding to the subgroup of $\pi_1 N$ generated by $\gamma$.
This covering space admits a canonical compactification $N_\gamma$, 
for example, using normal coordinates about the geodesic. Alternatively view $N_\gamma$ as the $\BZ$ quotient of $\BH^4\cup S^3_\infty$ by the loxodromic element corresponding to $\gamma$.  We identify $N_\gamma$ with $S^1 \times D^{3}$.  Since $f_k$ is homotopically trivial via a compactly supported homotopy, points of $\hat N$ are moved distances uniformly bounded above.  It follows that $\hat f_k$ extends to a homeomorphism $f^*_k$ such that $f_k^*|\partial N_\gamma=\id$.

Now the preimage of $\gamma$ in $\hat N$ consists of a single geodesic $\hat \gamma$ that maps 1-1 to $\gamma$ and infinitely many others that map $\infty$ to 1.  By \cite{BG} or using Proposition \ref{nullhfillings}, it follows that if $\gamma_i$ is one such lift, then $\hat f_k|N(\gamma_i)$ is isotopic to $\id$ via an isotopy that moves points uniformly bounded distance, independent of $i$.  Here $N(\gamma_i)$ is the corresponding lift of $N(\gamma)$.  The point here is that  $\delta_k$ is isotopically trivial when lifted to some finite sheeted cover.  Thus by an isotopy of $f_k^*$ which moves points uniformly bounded hyperbolic distance, we can assume that $f_k^*|N_\gamma$ is supported in $N(\hat\gamma)$, where we abuse notation by calling the isotoped map $f_k^*$.  I.e. $f^*_k$ is the standard  $\delta_k$ implantation in a neighborhood of the core geodesic of $N_\gamma$ . 

Finally  apply our $W_3'$-invariant to the resulting homeomorphism $f^*_k$ of $S^1 \times D^{3}$ to conclude that the $f_k$'s, $k\ge 4$, freely generate an infinite rank abelian subgroup  of $\pi_0(\Diff_0(N))$.  
\end{proof}
\end{thm}

The proof of Theorem \ref{mainthm} can be elaborated to construct a homomorphism $\pi_{n-4} \Homeo(N) \otimes \Rat$ to an
infinite direct sum of copies of $\Rat$ when $n>4$, in particular one copy of $\pi_{2n-3} C^\tau_3(S^1 \times D^{n-1}) \otimes \Rat / R'$
for each embedded orientable geodesic in $N$. The $n=4$ case is somewhat distinct, as the target group is a semi-direct product of
the group of hyperbolic isometries of $N$ and an infinite direct-sum of copies of $\Rat$. This follows from Mostow Rigidity and the finiteness
of the isometry groups of closed hyperbolic manifolds. The idea for the homomorphism is to find an infinite set $\gamma_i$ of distinct embedded geodesics and consider the diffeomorphisms $f_{i,k}$ obtained by implanting $\delta_k$ in $N(\gamma_i)$.  Again these diffeomorphisms generate a free abelian subgroup of $\pi_0(\Diff_0(N))$, provided $k\ge 4$.  The key point is that when lifting  to $\hat N_{\gamma_i}$ all the preimages of each $\gamma_j$, $j\neq i$ are non compact, so when extending to $N_{\gamma_i}$ only the $f_{j,k}$'s with $j=i$ survive up to isotopy and these are distinguished by $W'_3$.  I.e. our invariant of $\pi_{n-4} \Homeo(N) \otimes \Rat$ will
be $W_3'(f_{i,k})$, in the summand corresponding to $\gamma_i$.  

\providecommand{\bysame}{\leavevmode\hbox to3em{\hrulefill}\thinspace}

\end{document}